\documentclass[12pt]{article}

\title{Diophantine Approximation on varieties V: Algebraic independence criteria}

\author{Heinrich Massold}

\usepackage{bezier,amsfonts}

\usepackage{amsmath,amsfonts,bbm,amssymb}

\newtheorem{Satz}{}[section]

\newcommand{\satz}[1]{\vspace{2mm} \begin{Satz}{\bf #1}}

\newcommand{\proof}{\vspace{4mm} {\sc Proof}\hspace{0.2cm}}

\newcommand{\la}{\langle}

\newcommand{\ra}{\rangle}

\newcommand{\di}{{\mbox{div}}}

\newcommand{\R}{{\mathbbm R}}

\newcommand{\Pe}{{\mathbbm P}}

\newcommand{\Z}{\mathbbm{Z}}

\newcommand{\N}{\mathbbm{N}}

\newcommand{\vol}{\mbox{vol}}

\newcommand{\CO}{{\cal O}}

\newcommand{\CX}{{\cal X}}

\newcommand{\CY}{{\cal Y}}

\newcommand{\CZ}{{\cal Z}}

\newcommand{\CL}{{\cal L}}

\newcommand{\CF}{{\cal F}}

\newcommand{\fm}{{\mathfrak m}}

\newcommand{\Q}{{\mathbbm Q}}

\newcommand{\C}{{\mathbbm C}}

\newcommand{\A}{{\mathbbm A}}

\newcommand{\codim}{\mbox{codim}}

\newcommand{\spec}{\mbox{Spec}}

\begin{document}

\parindent0mm

\maketitle

\tableofcontents

\section{Introduction}

Let $\Pe^M$ be the projective space of dimension $M$ over $\spec \; \Z$,
and $\CX$ an irreducible arithmetic sub variety. A point $\theta \in \CX(\C)$ 
is called generic, if  the algebraic closure of $\{ \theta \}$
over $\Z$ is all of $\CX$. Part III of this
series of papers (\cite{App3}) established a lower bound for the 
approximability of generic points $\theta$ by algebraic points or sub varieties
in terms of the dimension of $\CX$,
which by \cite{Mahler} is best possible except for a subset 
of points of measure zero. More specifically, if the height and degree of 
an effective cycle on $\CX$ are
defined via $O(1)$ and $\overline{O(1)}$, and the algebraic distance 
of an effective cycle to $\theta$ is defined with respect to
$\mu = c_1(\overline{O(1)})$. (See \cite{App1}, Section4)

\satz{Theorem} \label{mainapp}
Let $\CX$ be an irreducible quasi projective arithmetic variety of relative 
dimension $t$ over $\spec \; \CO_k$, and $\bar{\CL}$ an ample positive 
metrized line bundle on some projective compactification of $\CX$. There is a 
number $b>0$ such that for every $a>>0$, and every generic 
$\theta \in X(\C_{\sigma})$ there is an infinite subset $M \subset \N$ such that 
for each $D \in M$ there is an irreducible subscheme $\alpha_D$ of codimension
$t$ fulfilling
\[ \deg \alpha_D \leq D^t, \quad h(\alpha_D) \leq a D^t, \quad
   \log |\alpha_D,\theta| \leq - b a D^{t+1}. \] 
%or otherwise
%\[ \omega^t(\theta) \leq \tilde{\omega}^t(\theta) \geq b. \]
\end{Satz}

\proof 
\cite{App3}, Theorem 1.2, Corollary 1.3.

\vspace{2mm}

It is the objective of this paper to reverse 
this conclusion, i.\@ e.\@ the approximability of a generic point by
algebraic subvarieties will imply a lower bound on the dimension of
$\CX$, and hence give criteria of algebraic independence of complex
numbers in terms of the approximability of corresponding points on arithmetic
varieties.

For the rest of the paper, if not specified otherwise,
$\partial_1, \ldots, \partial_t$ will be derivations 
of $k(X)$ whose restrictions to the
tangent space of $X$ at $\theta$ are linearily independent. For
a multi index $I = (i_1, \ldots, it) \in \N^t$, denote by $|I|$ its norm
$i_1+\cdots+i_t$ and by $\partial^I$ the differential operator
$\partial_1^{i_1} \cdots \partial_t^{i_t}$. Further, for a global
section $f \in \Gamma(X,O(D))$ denote
$|f|_{L^2(\Pe^M)} = (\int_{\Pe^M} |f|^2 \mu^M)^{1/2}$.

\satz{Theorem} \label{algind1}
Let $\CX$ be an irreducible subvariety of relative dimension $t$ in $\Pe^M$, and
$\theta = [(\theta_0,\ldots,\theta_m)] \in X(\C_\sigma)$ a generic point.
One may assume $\theta_0 \neq 0$, and then 
$t = \mbox{trdeg}_k (\theta_1/\theta_0, \ldots, \theta_M/\theta_0)$.
Let further, $D_k,S_k$ be series of natural numbers, $H_k,V_k$ series of
positive real numbers such that $S_k \leq D_k$, the series
$D_k/S_k, H_k/S_k, V_k/S_k$ are non- decreasing, and 
\[ \limsup_{k \to \infty} \frac{S_k^s V_k}{D_k^s (D_k + H_k)} = \infty. \]
Additionally assume that for each sufficiently big $k \in \N$, there is a set
of global sections $\CF_k$ of $O(D)$ such for each $f \in \CF_k$,
\[ \deg f \leq D_k, \quad \log |f|_{L^2(\Pe^>M)} \leq H_k, \quad
   sup_{|I| \leq S_k} \log |\partial^I (f/g^{\otimes D})(\theta)| \leq - V_k. \]
and that there is no point $x \in \Pe^M(\C)$ such that $f_x = 0$ 
for every $f \in \CF_k$, and $\log |x,\theta| \leq \frac{V_{k-1}}{S_{k-1}}$. 
Then $t$ is at least $s+1$.
\end{Satz}

The criterion entails the Philippon criterion if one takes $S_k=0$ for all $k$.
An alternative proof to the one given here was already given in \cite{LR} 
(Theorem 2.1). Under an additional assumption,
this new proof furthermore also entails a characterisation of the point
$\theta$ in terms of its approximability.

This criterion has a difficiency because it is usually used in cases in
which the series $(D_k,H_k,S_k,V_k)$ fulfill certain regularity conditions
(see below), and in this case there verifyably are points $\theta$
on any variety $\CX$ that fulfill the conclusion of Theorem \ref{algind1}
without fulfilling its premiss.

\satz{Definition}
A function $f: \N \to \N \; (\R)$ is said to be of 
uniform polynomial growth, if the limes
\[ n_f = \lim_{k\to \infty} \frac{k (f(k+1)-f(k))}{f(k)} \]
exists.
\end{Satz}

\satz{Lemma} \label{upg}

\begin{enumerate}

\item
The set of functions of uniform polynomial growth is closed under compositions,
sums, products, differences and quotients with 
\[ n_{f \circ g} = n_f n_g, \!\quad n_{f+g} = \mbox{max}(n_f,n_g), \!\quad 
   n_{fg} = n_f+n_g, \!\quad n_{1/f}=-n_f, \!\quad n_{-f} =n_f. \]
If $f$ is unbounded, $f^{-1}$ is defined via
\[ f^{-1}(n) = \inf \{ k | f(k) \geq n \}, \]
and $f$ is of uniform polynomial growth with $n_f \neq 0$,
then $f^{-1}$ is of uniform polynomial growth with $n_{f^{-1}} = 1/n_f$.

\item
If $f,g$ are of uniform polynomial growth, and $f(k) \geq g(k)$ for 
every sufficiently big $k$, then $n_f \geq n_g$.

\item
A function $f$ is of uniform polynomial growth, if and only if there is an
$n_f \in \R$ such that for every $\epsilon>0$, there is a $k_0 \in \N$ such 
that
\[ k^{n_f-\epsilon}< f(k) < k^{n_f+\epsilon} \]
for every $k \geq k_0$.

\item
If $f$ is of uniform polynomial growth, and $n$ is any natural number, then
for sufficiently big $k$,
\[ f(k+n) \leq 2 f(k). \]

\end{enumerate}
\end{Satz}

\satz{Definition} \label{reg}
Let $(D_k,S_k,H_k,V_k)$ be a quadrupel of sequences of natural
and positive real numbers with $S_k \leq D_k/3$. 
The quadrupel is said to be of regular polynomial growth if $D_k/S_k$ and 
$H_k/S_k$ are monotonously increasing and unbounded, and
the functions $f(k) = D_k/S_k$ and  $g(k)=H_k/D_k$ are of 
uniform polynomial growth with $n_f > 0$, and $g(k) \geq c >0$ for 
sufficiently big $k$.
\end{Satz}

\satz{Proposition} \label{S}
In the situation of Theorem \ref{algind1}, if additionally $(D_k,S_k,H_k,V_k)$
is of regular polynomial growth, and $trdeg_k(\theta) = s+1$, then
$\theta$ is an $S$-point in the sense of Mahler classification
\end{Satz}

\proof
\cite{Mahler}

%\vspace{2mm}
%For a subvariety $Y \subset \Pe^M$, let $I_Y(D) \subset \Gamma(\Pe^M,O(D))$ be
%the space of global section, whose restriction to $Y$ are $0$, and 
%for a global section $f \in \Gamma(\Pe^M,O(D))$, denote by $f_Y^\bot$ the
%orthogonal projection modulo $I_Y(D)$, and by $|f_Y^\bot|$ its norm.

\satz{Theorem} \label{algind2}
Let $\CX$ be a subvariety  of relative dimension $t$ of $\Pe^M$,
and $\theta \in \CX(\C)$ a generic point.     
Further, $D_k, H_k, S_k, V_k$ a quadrupel of sequences of natural and 
positive real numbers that is of regular polynomial growth, and fulfills
\[ \lim_{k \to \infty} \frac{S_k^s V_k}{D_k^s (D_k + H_k)} = \infty. \]
Assume that for every sufficiently big $k$, there is a set
$\CF_k \subset \Gamma(\Pe^M,O(D))$, such that
for every irreducible subvariety $\CY \subset \CX$ that has sufficiently
small distance to $\theta$, there is an
$f \in \CF_k$, and an $I$ with $|I| \leq S_k/3$ 
such that the restriction of $\partial^I f$ to $\CY$ is nonzero, and
\[ \log |f_k| \leq H_k , \quad 
   \sup_{|I| \leq S_k} \log |\partial^I (f/g^{\otimes D}(\theta)| \leq - V_k. \]
Then, $t \geq s+1$.

\end{Satz}

{\bf Remark:}
The conditions in 
Theorem \ref{algind2} are fulfilled e.\@ g.\@ if for every sufficiently
big $k$ there are $t$ global sections $f_1, \ldots, f_t$ of $\CL^{\otimes D_k}$
with
\[ \log |f_i| \leq H_k, \quad D^{S_k}(\di f_i,\theta) \leq - V_k, \quad
   i=1, \ldots, t, \]
and numbers $I_1, \ldots, I_t$ with $I_i \leq S_k$ such that the divisors of
the sections \\
$\partial^{I_1} f_1, \ldots, \partial^{I_t} f_t$ intersect properly. Another
important case, where the conditions of the Theorem are fulfilled, will
be when the global sections with small algebraic distance are obtained
by having high order of vanishing at a certain point, and behave well
with respect to differentiation.

\section{Prerequisites}

\satz{Lemma} \label{bashoeh}
Let $\CX$ be a regular projective arithmetic variety,
$\bar{\CL}$ a metrized line bundle on $\CX$, and $f$ a global section
of $\CL^{\otimes D}$. Then, for every effective cycle $\CZ$ on $\CX$ such that
the intersection of $\CZ$ with $\di f$ is proper,
\[ h(\di f . \CZ) = D h(\CX) + \int_Z \log |f| c_1(\bar{L})^m, \]
where $m$ is the dimension of $Z$. In particular, if
$\CZ$ is an effective cycle of pure codimension in projective space, and 
$f \in \Gamma(\Pe^t,O(D))$, then
\[ h(\di f) = D h(\CZ) + \int_Z \log |f| \mu^m, \]
with $\mu = c_1(\bar{L})$.
\end{Satz}

\proof
\cite{BGS}, Proposition 3.2.1.(iv).

\satz{Lemma} \label{funkhoeh}
Let $\CX, \CY$ be regular projective
arithmetic varieties, and $f: \CX \to \CY$ a morphism.
Then for every metrized line bundle $\bar{\CL}$ on $\CY$, and every 
cycle $\CZ$ in $\CX$, such that $\dim f (\CZ) = \dim \CZ$,
\[ h_{f^* \bar{\CL}} (\CZ) = h_{\CL} (f_* \CZ). \]
\end{Satz}

\proof
\cite{BGS}, Proposition 3.2.1.(iii).

\satz{Lemma} \label{projhoeh}
For $m< n$ let $\Pe^m \subset \Pe^n$ the projective subspace corresponding to 
a choice of $m+1$ homogeneous coordinates, and $\Pe^{n-m-1} \subset \Pe^n$
the subspace corresponding the remaining $n-m$ coordinates. With $\pi$
the map
\[ \pi : \Pe^n \setminus \Pe^m \to \Pe^{n-m-1}, \quad
   [v + w] \mapsto [v], \quad [v] \in \Pe^m, [w] \in \Pe^{n-m-1}, \]
and any cycle $\CZ$ in $\Pe^n$, such that $Z$ does not meet $\Pe^m$,
\[ h(\pi_* (Z)) \leq h(Z). \]
\end{Satz}

\proof
\cite{BGS}, (3.3.7).

\satz{Lemma} \label{Normrel}
For every $f \in \Gamma(\Pe^t_\C,O(D))$,
\[ \log |f|_\infty - \frac D2 \sum_{m=1}^t \frac 1m \leq
   \int_{\Pe^t_\C} \log |f| \mu^t \leq \log |f|_{L^2} \leq
   \log |f|_\infty. \]
\end{Satz}
\proof \cite{BGS}, (1.4.10).

\satz{Lemma} \label{polprod}
Let $f \in \Gamma(\Pe^t,O(D)), g \in \Gamma(\Pe^t,O(D'))$. Then,
\[ \log |f|_{L^2} + \log |g|_{L^2} -c_2(\log D+ \log D') \leq
   \log |f g|_{L^2} \leq \]
\[ \log |f|_{L^2} + \log |g|_{L^2}+ c_1(D+D') + \log {D+D'+t \choose t}. \]
\end{Satz}

\proof \cite{App2}, Lemma 3.2.

\vspace{2mm}

For $p \leq t$, $Z \in Z^p_{eff}(\Pe^t)$, and $\theta$ a point not contained
in the support of $Z$ in \cite{App1} the algebraic distance $D(Z,\theta)$,
is defined. Recall also the definition of the
derivated algebraic distance of a point $\theta$ to an effective
$X$ cycle in $\Pe^t$, whose support does not contain $\theta$ in \cite{App4}:
Let $I =(1_1, \ldots, i_{2t}) \in \N^{2t}$ denote a multi index, 
$|I| = i_1 + \cdots + i_{2n}$ its norm, and $\partial^I$ the differential 
operator $\partial^{i_1}/\partial x_1 \partial^{i_2}/\partial y_i \cdots
\partial^{i_{2t}}/\partial y_t$, and let $\varphi: \A^t(\C) \to \Pe^t(\C)$
be the affine chart with $\varphi(0) = \theta$. The derivated algebraic 
distance $D^S(Z,\theta)$ of order $S \in \N$ is defined as
\[ D^S(\theta,X) := sup_{|I|\leq S} \log |\partial^I \exp D(\theta,X)|. \]
If $\psi$ is another affine chart centered at $\theta$, the derivated algebraic
distance with respect to $\psi$ differs from that with respect
to $\varphi$ only by a constant depending on $\psi$ and $\varphi$ times
$S \log \deg X$. See \cite{App4}.

There are the following Propositions for the derivated algebraic distance.

\satz{Proposition} \label{glschn}
For $s,D \in \N$, and $f \in \Gamma(\Pe^t,O(D))$ let $F$ be the polynomial
of degree at most $D$ in $t$ variables that corresponds to $f$ with respect
to affine coordinates of $\Pe^t$ centered at $\theta$. Then,
with some positive constant $c$ only depending on $t$,
\[ D^S(\di f,\theta) = \sup_{s \leq S |J| = s} 
  \log \left| \left( \frac{\partial^s}{(\partial z_1)^{j_1} \cdots
  (\partial z_t)^{j_t}} F \right) (0) \right| - \log |f| + 
  O((S+D) \log (SD)). \]
\end{Satz}

\proof
\cite{App4}, Theorem 1.3.

\satz{Corollary}

\end{Satz}
In the situation of the Lemma,
\[ D^S(\di f,\theta) \leq \sup_{|J| \leq S} \log
   \left| (\partial^J F)(0)\right| + c(S+D) \log (SD). \]
\vspace{2mm}

\proof
Follows from the estimate
\[ \log |f| \geq - c D \]
for global sections $f$ of $O(D)$ with a fixed positive constant $c$.

We will need two special cases of the derivative metric
B\'ezout Theorem, proved in \cite{App4}, namely

\satz{Theorem} \label{DMBT}
Let $\CX,\CY$ be properly intersecting effective cycles in projective
space $\Pe^t_\Z$, and $S,\bar{S}$ natural numbers with
$S \leq \deg X/3, \bar{S} \leq \deg Y/3$. There is a positive constant
$d$ only depending on $t$, and a function $f$ from the set of natural numbers
less or equal $\deg X + \deg Y$ to the set of pairs of natural nubers less
or equal $\deg X$ and $\deg Y$ respectively, such that $pr_1 \circ f$ and
$pr_2 \circ f$ are surjective, and for every $\theta \in \Pe^t(\C)$ not
contained in the support of $X . Y$.

\begin{enumerate}

\item
For a given $k_0 \leq \deg Z_0 + \deg Z_1$, and
any $k \leq \deg Z_0 + \deg Z_1$ greater or equal $k_0$,
and $(\bar{\nu}_0,\bar{\nu}_1) = f(k)$,
\[ 2(\bar{\nu}_0-\nu_0) (\bar{\nu}_1-\nu_1) 
   \log |Z_0+Z_1,\theta| + 2D^S(Z_0.Z_1,\theta)+ 2 D(Z_0,Z_1) \leq \]
\[ (\bar{\nu}_0 - \nu_0) D^{3\nu_1}(Z_1,\theta) + 
   (\bar{\nu}_1-\nu_1) D^{3\nu_0}(Z_0,\theta) + \]
\[ O((\deg Z_0 \deg Z_1+ S) \log (S \deg Z_0 \deg Z_1)). \]

\item
\[ 2D(X,Y) + 2D^{S \bar{S}} (X . Y,\theta) \leq \]
\[ \mbox{max}(\bar{S} D^{3S}(X,\theta),S D^{3\bar{S}}(Y,\theta)) +
   d (\deg X \deg Y) \log (\deg X \deg Y), \]
and
\[ 2 D(X,Y) + 2 D(X. Y, \theta) \leq \]
\[ \mbox{max}(\bar{S} D(X,\theta), D^{3\bar{S}}(Y,\theta)) +
    d (\deg X \deg Y) \log (\deg X \deg Y). \]
\end{enumerate}
\end{Satz}

\proof
\cite{App4}, Theorem 1.9, Corollary 1.11.

\satz{Corollary} \label{DMBTcor}
\begin{enumerate}

\item

For $S_0, d_0 \leq \deg Z_0/3,S_1 \leq Z_1/3$, and $S=S_0 S_1$, there is a 
$K \leq d_0 S_1$ such that
\[ K \log |Z_0+Z_1,\theta| + 2D^S (Z_0.Z_1,\theta)+ 2 D(Z_0,Z_1) \leq \]
\[ \mbox{max}(S_1 D^{9S_0}(Z_0,\theta) , d_0 D^{9S_1}(Z_1,\theta) + \]
\[ O((\deg Z_0 \deg Z_1+ S) \log (S \deg Z_0 \deg Z_1)). \]

\item
For $S_0 \leq \deg Z_0/3$, $S_1 \leq \deg Z_1/3$, and 
$|Z_0,\theta| \leq |Z_1,\theta|$,
\[ 2 D^{S_1}(Z_0,Z_1) + 2 D(Z_0,Z_1) \leq \]
\[ D^{3S_1}(Z_1,\theta) + 
   O((\deg Z_0 \deg Z_1+ S) \log (S \deg Z_0 \deg Z_1)). \]

\end{enumerate}
\end{Satz}

\proof
The proof is similar to the one of \cite{App4}, Corollary 1.11.

\vspace{2mm}

Similar to \cite{App3}, Proposition 2.4.1, on can also deduce

\satz{Theorem} \label{DMBText}
Let $\CY$ be an irreducible effective cycle of codimension $p$
in projective space, 
$f \in \Gamma(\Pe^t,O(D))_\Z$ a global section whose restriction to
$\CY$ is nonzero, and $\bar{f} \in \Gamma(\Pe^t,O(D))_\R$ a global section
that is orthogonal to $I_\CY(D)$ the elements of degree $D$ in the 
vanishing ideal of $Y$ such that $f_Y = \bar{f}_Y$. Then for  natural 
numbers $S\leq \deg Y/3,\bar{S} \leq D/3$ such that for every 
$\theta \in \Pe^t(\C)$ not contained in $\di f . Y$,
\[ 2 D^{S \bar{S}} (Y . \di f, \theta) \leq \]
\[ \mbox{max}(\bar{S} D^S(\di \bar{f} ,\theta), S D^{\bar{S}}(Y,\theta)) +
   \deg Y \log |f^\bot_Y| + D h(\CY) + 
   d D \deg Y \log (D \deg Y), \]
and
\[ D(Y . \di f, \theta) \leq \mbox{max}(\bar{S} D(\di \bar{f}, \theta),
   D^{\bar{S}} (Y,\theta)) + d D \deg Y \log(d D \log Y). \]
\end{Satz}

%%The concept of the drivated algebraic distance will now be slightly
%extended. Let $\Pe^M_\C$ be projective space of dimension $M$, and
%$X \subset \Pe^M\C$ a subvariety of dimension $t$. The canonical
%K\"ahler structure on $\Pe^M$ induces a K\"ahler structure on $X$; for
%every $\theta \in X$, and $U_\theta$ a neighbourhood of $\theta$ in $X$
%denote by $D_1^S U_\theta$ the space of parellel differential forms
%of order $S$ on $U_\theta$.

%\satz{Proposition}
%Let $X$ be a subvariety of dimension $t$ of $\Pe^M$, and $Y$ an effective
%cycle of pure codimension $p$ on $X$. Further let $\theta \in X$ be a 
%point not contained in the support of $Y$, and denote by $\Pe_X$ the
%tangent space of $X$ at $\theta$

%\begin{enumerate}

%\item
%For every subspace $\Pe(F) \subset \Pe_X$ properly intersecting

%\end{enumerate}

%\end{Satz}

Let $H: \N \to \R$ be function of uniform polynomial growth such that
$H(D)/D \geq a>0$ with a sufficiently big constant $a$, hence by Lemma
\ref{upg}.4, $n_H \geq 0$.

For $\CX$ an effective cycle in $\Pe^t$ define the $H/D$-size of $\CX$ as
\[ t_{\frac HD}(\CX) = \frac HD \deg X + h(\CX). \]

%\satz{Proposition} \label{1hilf}
%There are constants $c,b>0$ only depending on $t$ such that if
%\[ V \geq c \frac{D^t (H+D)}{S^t}, \]
%and there is a $f \in \Gamma(\Pe^t,O(D))$ such that
%\[ \log |f| \leq H, \quad D(\di f, \theta) \leq - V, \]
%then there is an irreducible zero dimensional cycle $\alpha$ in $\Pe^t$
%such that 
%\[ \deg \alpha \leq (D/S)^t, \quad h(\alpha) \leq (D/S)^{t-s}(H/S); \quad
%   \varphi_{\theta,H}(\alpha) \leq -b t_H(\alpha) \frac DS 
%   \frac{V S^t}{D^t(D+H)}. \]
%\end{Satz}

\satz{Proposition} \label{2hilf}
There are constants $c,b,\bar{b}>0, n \in \N$ only depending on $t$ such that 
for every generic $\theta \in \Pe^t(\C)$ and every function $H:\N\to \R$ as
above, there is an infinite 
set $M \subset \N$ such that for every $D \in \N$, there is an irreducible
zero dimensional subvariety $\alpha_{nD}$ of $\Pe^t_\Z$, a locally complete
intersection $\CX$ of codimension $s\leq t-1$ at $\alpha_D$ and global sections
$f \in \Gamma(\Pe^t,O(D))_\Z, \bar{f} \in \Gamma(\Pe^t,O(D))_\R$
such that $f^\bot_{\alpha_{nD}} = \bar{f}^\bot_{\alpha_{nD}} \neq 0$, and
with $\CX_{min}$ the irreducible component of $\CX$ with minimal $H/D$-size,
\[ \deg \CX \leq D^s, \quad h(\CX) \leq HD^{s-1}, \]
\[ \log |\bar{f}_{\alpha_{nD}}^\bot| \leq H, \quad
   \log |\la f | \theta \ra| \leq - b t_{\frac HD}(\CX_{min}) D^{t+1-s}, \]
\[ \deg \alpha_{nD} \leq (n D)^t, \quad h(\alpha_{nD}) \leq (nH) (nD)^{t-1}, 
   \quad D(\alpha_{nD},\theta) \leq - \bar{b} t_{\frac HD}(\alpha_{nD}) D, \]
\[ \log |\alpha_{nD},\theta| \leq - \bar{b} t_{\frac HD}(\alpha_{nD}) D, \quad
   t_{\frac HD}(\alpha_{nD}) \geq c t_{\frac HD}(\CX_{min}) D^{t-s}. \]
\end{Satz}

\proof
\cite{App3}, Corollary 4.21. One has to be cautious to adjust the constants.

%\satz{Lemma} \label{Exfs}
%Let $\CY$ be a successive $(D,a)$-intersection of
%codimension $s \leq t-1$ in $\Pe^t$. There are $n \in \N, b>0$ only depending
%on $t$, and a locally complete $D$-intersection $\CX$ at $\CY$ of codimension 
%$r \leq s$, and 
%global sections $g \in \Gamma(\Pe^t,O(n D))_\Z$,
%$\bar{g} \in \Gamma(\Pe^t,O(n D))$, such that
%$g_Y^\bot = \bar{g}_Y^\bot\neq 0$, and 
%\[ \log |g_Y^\bot| \leq 6 a n D, \quad
%   \log |\la \bar{g}|\theta \ra| \leq - 
%   b t_a(\CX_{min}) D^{\dim \CX}, \]
%and the restrictions of $f$, and $\bar{f}$ to every irreducible
%component of $X$ are non zero.
%\end{Satz}

%\satz{Lemma} \label{zwischen}
%For every irreducible $(D,a)$-approximation cycle $\CY$ of codimension
%$s \leq t-1$, belonging to an approximation chain $\CC$, either there is a 
%$(D,a)$-approxi- mation chain
%\[ \CC_1: \quad \Pe^t = \CY_0 \supset \cdots \supset \CX, \]
%such that $\CX$ is a locally complete $D$-intersection at $\CY$, and
%$\CC_1 \prec \CC$, or there is a $(D,a)$-intersection chain
%\[ \CC_2: \quad \CY_0 \supset \cdots \supset \CY \supset \CZ \]
%with $\CC_2 \prec \CC_1$.
%\end{Satz}

%\proof
%\cite{App3}, Lemma 

\vspace{2mm}

Another important tool for the proofs is the Liouville inequality.

\satz{Proposition: Liouville inequality} \label{liou}
Let $f \in \Gamma(\Pe^t,O(D))$, and $\alpha$ an algebraic point such
that $f_\alpha \neq 0$. There is a constant $d$, only depending on $t$
such that
\[ D(\di f, \alpha) \geq 
   D h(\alpha) - \deg \alpha \log |f| - d D \deg \alpha. \]
\end{Satz}

For the relation of this Proposition to the classical formulation of
the Liouville inequality, compare \cite{Liouville}.

\proof
Since by \cite{App1}, Theorem 2.2.2, $h(\di f) \leq \log |f| + D \sigma_t$,
this is a special case of the equality
\[ D(\di f, \alpha) = h(\di f . \alpha) - \deg \alpha h(\di f) -
   \deg f h(\alpha) + \sigma_t \deg f \deg \alpha \]
from \cite{App1}, Scholie 4.3, together with the estimate
$D(\di f, \alpha) \leq d' \deg \alpha \deg f$ from \cite{BGS}, Proposition
5.1.

\section{Derivatives}

\subsection{Polynomials modelling derivatives of rational functions}

With $\CX$ an arithmetic sub variety of relative dimension $t$ in $\Pe^M_\Z$,
and $g$ a global section of $O(1)$ whose restriction to $\CX$ is nonzero, 
let $\theta \in X(\C_\sigma)$ be a generic point, and 
$\partial_1, \ldots, partial_t$ derivatives of $X$ as in the introduction.

\satz{Lemma} \label{dadrf}
With the above notations, and $f$ a global section of $\CL^{\otimes D}$,
\[ \sup_{|I| \leq S} \log 
   \left|\partial^I \frac{f}{g^{\otimes D}}(\theta) \right| =
   D^S(\di f, \theta) + \log |f|_{L^2(\Pe^M)} + O((S+D) \log SD), \]
for every $S \leq D$.
\end{Satz}

\proof
Let $U_\theta$ be a neighbourhood of $\theta$, and
$\varphi: U \to U_\theta$ an affine chart of $U_\theta$. Further,
$\tilde{\partial}_1, \ldots \tilde{\partial}_t$ the canonical derivatives
on $U$. Then, 
\[ (\varphi^{-1})^* \partial = h \tilde{\partial} \]
with a $(t \times t)$-matrix of rational functions $h$. Hence,
\[ \sup_{|I| \leq S} \log 
   \left|\partial^I \frac{f}{g^{\otimes D}}(\theta) \right| =
   \sup_{|I| \leq S} \log
   \left|\tilde{\partial}^I (\varphi^* f) (0) \right| +
   O((S+D) \log SD), \]
and the Lemma follows from \cite{App4}, Theorem 1.3.  

\satz{Corollary}
If $\bar{g}$ is another global section of $\CL$, and
$\bar{\partial}_1, \ldots, \bar{\partial}_t$ another set of derivations
of $k(X)$ whose restrictions to $T_\theta X$ are linearily independent, then
\[  \sup_{|I| \leq S} \log 
   \left|\partial^I \frac{f}{g^{\otimes D}}(\theta) \right| =
   \sup_{|I| \leq S} \log 
   \left|\bar{\partial}^I \frac{f}{\bar{g}^{\otimes D}}(\theta) \right| 
   + O((S+D) \log SD). \]
\end{Satz}

\vspace{2mm}

Because of this Corollary to the derivatives of a global section it doesn't
matter which derivatives $\partial_1,\ldots,\partial_t$ in $k(X)$ one choses.
In the proofs of the main Theorem we will chose them according to the
definition in the next Theorem.

Let $\CX \subset \Pe^M_\Z$ be 
an irreducible subvariety of relative dimension $t$, and $\Pe^t \subset \Pe^M$
a subspace defined over $\spec \Z$ such $(\Pe^t)^\bot$ does not meet $\CX$.
Then, with
\[ \pi : \Pe^M \setminus (\Pe^t)^\bot \to \Pe^t, \]
the restriction $\pi_X$ of $\pi$ to $\CX$ is a proper map from $\CX$ to 
$\Pe^t$. Denote by $x_0, \ldots, x_M$
the homogeneous coordinates of $\Pe^M$ ordered such that
$x_0,\ldots,x_t$ are homogeneous coordinates of $\Pe^t$.
There is the canonical map $k(x_1, \ldots, x_M) \cong k(\Pe^M) \to k(X)$, and
if $\bar{x}_i, i=1, \ldots, t$ denotes the image of $x_i$ under this map,
the function field $k(X)$ is a finite extension of 
$k(\bar{x}_1, \ldots, \bar{x}_t)$.

We denote by $\partial/\partial x_\mu$ the usual derivations 
of $k(x_1, \ldots,x_M) \cong k(x_1/x_0,\ldots,x_M/x_0)$, and do not distinguish 
between a polynomial
$f(x_1, \ldots, x_M)$, its image \\ $f(x_1/x_0, \ldots,x_M/x_0)$ in $k(\Pe^M)$, 
and its image $f(\bar{x}_1, \ldots, \bar{x}_M)$ in $k(X)$.
The following Theorem is a generalization of \cite{RW}, Proposition ??  to
higher dimensions.

\satz{Theorem} \label{Abl}
With the above notations,
let $\partial_t, \ldots, \partial_t$ be the derivations of $k(X)$ defined by
\[ \partial_l x_l =1, \quad \mbox{and} \quad
   \partial_l x_i = 0 \quad \mbox{for} \quad 
   i \in \{ 1, \ldots,t \} \setminus \{l\}. \]
Let further $I =(i_1, \ldots, i_t) \in \N^t$ be a multi index of degree
$S=i_1+\cdots+i_t$, and 
$\partial^I = \partial_1^{i_1} \cdots \partial_t^{i_t}$.

There is a homogeneous polynomial $P=P(x_0, \ldots x_M)$ with
\[ \deg P \leq (M-t) \deg X, \quad
   \log |P|_{L^2(\Pe^M)} \leq c \deg X + h(\CX), \]
with $c$ a constant only depending on $M$ and $t$, such that for every 
multi index $I$ of degree $S$, and every homogeneous polynomial $f$, 
\[ \partial^I f = \frac{f_I}{P^{2S-1}}. \]
where $f_I$ is a homogeneous polynomial with
\[ \deg f_I \leq \deg f + (2S-1)(M-t) \deg X, \]
\[ \log |f_I|_{L^2(\Pe^M)} \leq \]
\[ \log |f| + \log \deg f+
   (2S-1)(M-t) (h(\CX) + c \deg X +\log \deg X) + \log (2S!). \]
\end{Satz}

\proof
Let $\pi_X$ be the restriction of $\pi$ to $\CX$. For 
any $\mu = t+1, \ldots, M$, the projection of $\CX$ to
the space with homogeneous coordinates $x_0, \ldots, x_t, x_\mu$ is a 
hyper surface of
degree at most $\deg X$. Let $P_\mu$ be the corresponding homogeneous 
polynomial in $x_0, \ldots, x_t, x_\mu$. Then
$\deg P_\mu \leq \deg X$, and by Lemma \ref{bashoeh}, \ref{projhoeh} and
\ref{Normrel},
\begin{equation} \label{Pabsch}
\log |P_\mu| \leq \int_{\Pe^M} \mu^M + c \deg X 
\leq h(\pi_* \CX) + c \deg X \leq h(\CX)+ c \deg X. 
\end{equation}
Let further 
\[ A_0 := \prod_{\mu=t+1}^M \frac{\partial P_\mu}{\partial x_\mu}, \]
and
\begin{equation} \label{abl0}
 A_{l \mu} := - \frac{\partial P_\mu}{\partial x_l}
   \left(\frac{\partial P_\mu}{\partial x_\mu}\right)^{-1} \in k(X), 
\end{equation}
for $l=1, \ldots,t$, and $\mu=t+1, \ldots, M$. 
We have 
\[\deg A_0 \leq (M-t)(\deg X-1), \]
and using Lemma \ref{polprod}, and (\ref{Pabsch}),
\[ \log |A_0| \leq (M-t) (h(\CX) +c \deg X + \log \deg X). \] 
Also, $A_0 A_{l \mu}$ is
a polynomial with 
\[ \deg (A_0 A_{k \mu}) \leq (M-t) (\deg X -1), \]
\[ \log |A_0 A_{k \mu}| \leq (M-t) (h(\CX) + 
   c \deg X +\log \deg X). \] 
Since $P_\mu(x_1, \ldots, x_t,x_\mu) =0$ on $X$, we get
\[ 0 = \partial_l P_\mu =
   \frac{\partial P_\mu}{\partial x_l} + 
   \frac{\partial P_\mu}{\partial x_\mu} \partial_l x_\mu, \]
hence,
\[ \partial_l x_\mu = - \frac{\partial P_\mu}{\partial x_l}
   \left( \frac{\partial P_\mu}{\partial x_\mu} \right)^{-1} =
   A_{l \mu}, \]
and
\[ f_l = A_0 \frac{\partial f}{\partial x_l} +
   \sum_{\mu=t+1}^M A_0 A_{l \mu} \frac{\partial f}{\partial x_\mu} \]
is a polynomial with
\[ \deg f_l \leq \mbox{max}(\deg A_0,\deg (A_0 A_{l\mu})) +\deg f-1 \leq
   (M-t)\deg X + \deg f, \]
\begin{equation} \label{abl1}
\log |f_l| \leq 
\log \deg f + \log |f| + (M-t)(h(\CX) + c \deg X +\log \deg X) + 
\log 2,
\end{equation}
and
\begin{equation} \label{abl2}
\partial_l f = \frac{f_l}{A_0}.
\end{equation}
Put $P = A_0$. Then $\deg P \leq (M-t) (\deg X-1)$, and the estimate on the 
norm of $P$ immediately follows from (\ref{Pabsch}), and Lemma \ref{polprod}.

Assume now the Theorem proved for $I$ of degree $S$. 
That is, for any $I$ with $|I|=S$,
\[ \partial^I f = \frac{f_I}{P^S}, \]
for some polynomial $f_I$ with norm and degree
fulfilling the estimates from the Theorem.
Then, with $\bar{I} = I + (0, \ldots, 0,1,0, \ldots,0)$,
\[ \partial^{\bar{I}} f=\partial_l \partial^I f = \partial_l 
   \frac{f_I}{P^{2S-1}} = 
   \frac{(\partial_l f_I) P^{2S-1} - f_I (2S-1) P^{2S-2} \partial_l P}
   {P^{4S-2}} = \]
\[ \frac{P^2 \partial_l f_I - (2S-1)f_I P \partial_l P}{P^{2S+1}}. \] 
By (\ref{abl1}), (\ref{abl2}) and induction hypothesis 
$P^2 \partial_l f_I$, and $(S-1) P f_I \partial_l P$ are polynomials with
\begin{eqnarray*} 
\deg(P^2 \partial_l f_I) \leq \deg P + \deg f_I &\leq&
   2 (M-t) \deg X + (2S-1)(M-t) \deg X \\
&+& \deg f \\
   &=& (2 S+1)(M-t) \deg X + \deg f, 
\end{eqnarray*}
and 
\begin{eqnarray*} \deg (f_I P \partial_l P) &\leq& 
   2(M-t) \deg X + (2S-1)(M-t) \deg X + \deg f \\
&=& (2S+1)(M-t) \deg X + \deg f. 
\end{eqnarray*}
Likewise, the norms of
$P^2 \partial_l F_I$, and $(S-1) P f_I \partial_l P$ by (\ref{abl1});
(\ref{abl2}) and induction hypothesis fulfill the inequalities
\begin{eqnarray*} 
\log|P^2 \partial_l f_I| &\leq& (2S+1)(M-t) (h(\CX) +
   c \deg X +\log \deg X) +\log (2S)! \\
&+& \log |f| + \log \deg f \\ \\
&\leq&2(S+1) (M-t) (h(\CX) + c \deg X + \log \deg X) + 2S \log 2S^2 \\
&+&  \log |f| + \log \deg f, 
\end{eqnarray*}
and
\begin{eqnarray*}
\log |(2S-1)f_I P \partial_k P| & \leq &
   \log (2S-1) + (2S-1)(M-t) \times \\
&&(h(\CX) + c \deg X +\log \deg X) \\
&+& \log |f| + \log \deg f + \log (2S)!.
\end{eqnarray*}
Hence, with $f_{\bar{I}} = P^2 \partial_l f_I - (2S-1)f_I P \partial_l P$, we
have $\partial^{\bar{I}} f = f_{\bar{I}}/P^{2S+1}$, and
\[ \deg f_{\bar{I}} \leq (S+1)(M-t) \deg X + \deg f, \]
\[ \log |f_{\bar{I}}| \leq  (2S+1) (M-t) 
   (h(\CX) + c \deg X + \log \deg X) + \log (2S+2)!, \]
and the claim follows for $S+1$.

\satz{Corollary} \label{distAbl}
With the notations of the Theorem, for every $\theta$ in $\CX(\C)$, such
that $f(\theta) \neq 0$ and $P(\theta) \neq 0)$, there is a constant $c$ only 
depending on $\theta$, and $\CX$ such that 
\[ D^S(\di f, \theta) = \sup_{|I| \leq S} \log |f_I(\theta)| +
   O((S+D) \log(SD)). \] 
Moreover, for every $T \leq S$,
\[ D^S(\di f, \theta) = \sup_{|I| \leq S-T} \sup_{|J| \leq T}
   \log |(\partial^J f_I)(\theta)|. \]
\end{Satz}

\proof
Since $\log |P(\theta)^{2S-1}| = c (2S-1)$, for some constant $c$, 
with $g = x_0^D$, the claim follows from the Theorem, together with
Lemma \ref{dadrf}.

%Set $c = -\log |P(\theta)|$. There is a neighbourhood $U$ of $\pi \theta$ 
%such that $\pi_X^{-1}: U \to X(\C)$ is an affine chart.
%Hence, With $\tilde{\partial}^I$ the differential operator on $U$,
%\[ D^S(\di f , \theta) = \sup_{|I| = S} 
%   \log |(\tilde{\partial}^I (\pi_X^{-1})^*f)(\pi \theta)| = \sup_{|I|=S}
%   \log |((\pi_X)_*(\tilde{\partial}^I) f)(\theta)| = \]
%\[ \sup_{|I| \leq S} \log |(\partial^I f)(\theta)| =
%   \sup_{|I| \leq S} \log \left|\frac{f_I(\theta)}{P(\theta)^S} \right| =
%   c(2S-1) + \sup_{|I| \leq S} \log |f_I(\theta)|. \]

%\satz{Propostion} \label{gshr}
%With the notations of the Theorem, for every $f \in k(x_1, \ldots, x_t)$,
%\[ \partial^I (\pi^*_X) f = \frac{\partial^{|I|}}{(\partial x)^I} f, \]
%and for every $f \in k(x_1, \ldots, x_M)$,
%\[ \left|\frac{\partial^{|I|}}{(\partial x)^I} (\pi_X)_* f 
%   (\pi_X \theta)\right| \leq
%   \left| (\partial^I f) \theta \right|. \] 
%\end{Satz}

\subsection{Local B\'ezout Theorem}

In this subsection $k$ is a field of characteristic zero and $X$ a scheme of 
dimension $t$ over $\spec \; k$. 
For $y$ a point in $X$ denote by $Y = \overline{\{ y \}}$ its Zariski closure.

\satz{Definition}
\begin{enumerate}

\item
Let $y$ be a point in $X$ with
with $\dim Y =t-p$ .  For $Z$ an irreducible subvariety of codimension $p-1$,
$f \in k(Z)$ and $\mathfrak{m}_y \subset \CO_{X,y}$ the maximal ideal
in the localization of $\CO_X$ at $y$, 
define the order of vanishing $v_y(f)$ of $f$ at $\CY$ as
\[ v_y(f) := \mbox{max} \{ n \in \N| f \in \mathfrak{m}_y^n \}. \]

\item
For $X$ an irreducible subscheme of $\Pe^M_k$, and
$\Pe(W) \subset \Pe^t = \mbox{Proj} k[x_0, \ldots, x_M]$ a 
projective subspace of codimension $q$, let $w$ be the corresponding point 
in $\Pe^M$ and $Y$ an irreducible subvariety of 
codimension $p$ with $p \leq q$, define $v_w(Y)$ as
\[ v_w(Y) := \mbox{min}_{\Pe(F)} \{ \mbox{multiplicity of} \; \Pe(W) \;
   \mbox{in} \; \Pe(F) . Y\}, \]
where $\Pe(F)$ runs over all subspaces $\Pe(F) \subset \Pe^M$ of 
codimension $q-p$ that intersect $Y$ properly and contain $\Pe(W)$.
Define $v_w:Z(\Pe^t) \to \Z$ by linear extension.

\item
For $X = \Pe^M$, $w \in \Pe^M$ a point corresponding to a subspace 
$\Pe(W) \subset \Pe^M$, and $ZZ_1-Z_2$ a cycle of pure codimension $p$
in $\Pe^M$ define the order of vanishing of $Z$ at $w$ as the difference of
the orders of
vanishing as defined in part 1 of the chow forms $f_{Z_1},f_{Z_2}$ of 
$Z_1,Z_2$ at the subvariety
\[ \Pe(\check{W})_i = \check{\Pe}^M \times \cdots \times
   \check{\Pe}^M \times \Pe(\check{W}) \times \check{\Pe}^M \times \cdots
   \times \check{\Pe}^M, \]
where $\check{\Pe}^M$ is the space dual to $\Pe^M$, and $\Pe(\check{W})$
the space dual to $\Pe(W)$.
Since the chow divisor is invariant under permutation of the factors
in $(\check{\Pe}^t)^{M+1-p}$, this number does not depend on the choice
of $i \in \{1, \ldots, M+1-p \}$.
\end{enumerate}
\end{Satz}

\satz{Lemma}
\begin{enumerate}

\item
For $y_w$ the point corresponding to a subspace
$\Pe(W) \subset \Pe^M$ of codimension $q$, and $Z$ a subvariety of 
codimension $q-11$ in $\Pe^t$ the definitions in 1 and 2 coincide.

\item
The Definitions 2 and 3 coincide.

\end{enumerate}
\end{Satz}

%\proof

\satz{Fact} \label{dim0}
\begin{enumerate}
\item
If $w$ is a point corresponding to a subspace, $X$ an effective cycle in 
$\Pe^M$, then $\Pe(W) \subset \mbox{supp} \;X$, if and only if $v_w(X) \geq 1$.

\item
If $y$ is a closed point of $\Pe^M$, and $X$ an effective cycle of
pure codimension $M$, the multiplicity of $y$ in $X$ equals $v_y(X)$.

\item
Let $q \geq p$, and $\Pe(W),\Pe(F)$ be subspaces of codimension $q$, and
$p$ respectively. If $w$ is the point corresponding to $\Pe(W)$, then
\[ v_w(\Pe(F)) = 1 \Leftrightarrow \Pe(W) \subset \Pe(F), \quad \mbox{and}
   \quad v_w(\Pe(F)) = 0 \Leftrightarrow \Pe(W) \not\subset \Pe(F). \]

\item
Let $\Pe(W) \subset \Pe(F) \subset \Pe^M$ be subspaces, and
$Y$ an effective cycle intersecting $\Pe(F)$ properly. If
$v^{\Pe(F)}_w(Y)$ is defined as the order of vanishing of $\Pe(F) . Y$ at
$\Pe(W)$ inside $\Pe(F)$, then
\[ v^{\Pe(F)}_w(Y) \geq v_w(Y). \]

\item
Let $X$ be an effective cycle of pure codimension $p$ in $\Pe^t$, and $y$ a 
closed point. Then for every subspace $\Pe(F) \subset \Pe^t$, of 
codimension $t-p$ containing $y$, and intersecting $Y$ properly,
\[ v_y(X) \leq  v_y(\Pe(F). X), \]
and there exists a subspace $\Pe(F)$ with these properties such that
\[ v_y(X) = v_y(\Pe(F).X). \]

\end{enumerate}
\end{Satz}

\satz{Proposition} \label{vfunc}
Let $X$ be an irreducible subvariety of dimension $t$ of $\Pe^M$,
and $w$ a closed point in $X$. Further, $f,g \in \Gamma(\Pe^M,O(D))$ with 
$f_y \neq 0$. If for a natural number
$S$ and every multi index $I$ with $|I| \leq S$ the equality
$(\partial^I f) (y) = 0$ holds, then the order of vanishing
$v_y(Z)$ of $Z=X. \di f$ at $y$ is at least $S$.
\end{Satz}

\proof
By Fact \ref{dim0}, there is a subspace $\Pe(F) \subset \Pe^M$ be of 
codimension $t-1$ containing $y$ and properly intersecting $Z$
such that $v_y(Z) = v_y(\Pe(F).Z)$. Since $g_y \neq 0$,
the multiplicity of $y$ in $\Pe(F). X. \di f$ equals the multiplicity of
$y$ in $\Pe(F).X. \di(f/g)$, that is
\[ v_y(Z) = v_y(\Pe(F).Z) = v_y(\Pe(F).X.\di (f/g)). \]
Further, if $\bar{f},\bar{g}$ are the restrictions of $f,g$ to one-dimensional
subvariety $\Pe(F) \cap X$, then
\[ v_y(\Pe(F).X.\di (f/g)) \geq v_y(\di (\bar{f}/\bar{g})). \]
If $\partial$ is a derivation of $\Pe(F) \cap X$ whose restriction to
$y \in \Pe(F) \cap X$ is nonzero, then $\partial$ is a linear combination 
with coefficients in $k(X)$ of $\partial_1, \ldots, \partial_t$, hence
$\partial^s (\bar{f}/\bar{g}) = 0$ for every $s \leq S$, which means that
$(\bar{f}/\bar{g})$ is contained in the $S$th power $\fm_{\Pe(F)\cap X,y}^S$
of the maximal ideal $\fm_{\Pe(F)\cap X,y} \subset \CO_{\Pe(F)\cap X,y}$, that is
$v_y(\bar{f}/\bar{g}) \geq S$. Together with the above equalities and
estimates this implies the claim.

\vspace{2mm}

Two effective cycles $Y,Z$ of projective space are said to intersect
properly at a point $x \in \Pe^M$ if for every irreducible component
$U$ of the intersection of the supports of $Y$ and $Z$ that contains $x$
the equality $\codim W = \codim Y + \codim Z$ holds.

\satz{Local B\'ezout Theorem} \label{lbez}
For $x$ a closed point in $\Pe^M$ and two cycles $Y,Z$ of $\Pe^M$, intersecting
properly at $x$,
\[ v_x(Y.Z) \geq v_x(Y) v_x(Z). \]
\end{Satz}

\satz{Remark}  \label{luraum} 
By Fact \ref{dim0},
the Theorem holds in case $Y$ is a projective subspace $\Pe(F) \subset \Pe^t$.
\end{Satz}

\satz{Lemma}
Let $Y,Z$ be properly intersecting irreducible varieties of codimension $p,q$
of $\Pe^M$, and $X \# Y$ their join. For a closed point $x$ in $\Pe^m$ 
there are subpaces $\Pe(F), \Pe(F')$ of codimensions $t-p,t-q$ containing
$x$ such that the intersections $\Pe(F) .Y$ and $\Pe(F').Z$ are proper, and
\[ v_{(x,x)}(Y \# Z) = v_{(x,x)}(Y \# Z . \Pe(F) \# \Pe(F')). \]
\end{Satz}

\satz{Lemma}
A point $y\#z$ in $\Pe^{2M+1}$ intersects $\Pe(\Delta)$ if and only if
$y=z$. Further, $(y \# y) . \Pe(\Delta) = (y,y)$.
\end{Satz}

\proof
Let $u \in k^{t+1}, v \in k^{t+1}$ be vectors representing $y,z$, i.\@ e.\@
$[u] = y, [v]=z$. The join $y\#z \subset \Pe^{2t+1}$ consists of the points
$[(au,bv)]$, $a,b \in k$. If $[(au,bv)] \in \Pe(\Delta)$, then there is
a vector $w \in k^{t+1}$ such that $(au,bv) = (w,w)$. Hence,
$au = w=bv$, that is $y= [u] = [v] = z$, and $[(w,w)] = (y,y)$.

\satz{Lemma}
Let $x$ be a closed point in projective space, $Y,Z$ properly intersecting 
effective cycles in $\Pe^M$, and $Y \# Z$ their join in $\Pe^{2M+1}$.
\begin{enumerate}
\item
\[ v_{(x,x)}(Y \# Z) \geq v_x(Y) v_x(Z). \]
\item
\[ v_{(x,x)}(\delta_*(Y.Z)) = v_x (Y. Z), \]
where $\delta: \Pe^M \times \Pe^M \to \Pe(\Delta)$ is the diagonal embedding.
\end{enumerate}
\end{Satz}

\proof
1. By Fact \ref{dim0}.4, there are subspaces
$\Pe(F), \Pe(F') \subset \Pe^t$ such that
$v_x(Y) = v_x(\Pe(F).Y), v_x(Z) = v_x(\Pe(F').Z)$. 
Since $\Pe(F).Y = \sum_y n_y y$ is zero dimensional, by Lemma \ref{dim0},
$v_x(\Pe(F).Y) = n_x$ similarly, with $\Pe(F').Z = \sum_z m_z z$, the
equality $v_x(\Pe(F').Z) =  m_x$ holds.
Since $(Y \# Z).(\Pe(F) \# \Pe(F')) = \sum_{y,z} n_y n_z y \# z$, and
$x \# x$ contains $(x,x)$, it follows from the previous Lemma  that
\[ v_{(x,x)}(Y \# Z) = v_{(x,x)} (Y \# Z).(\Pe(F) \# \Pe(F')) \geq
   n_x m_x = v_x(Y) v_x(Z). \]

\vspace{2mm}

2. Since the diagonal embedding is an isomorphism, this follows from the
previous Lemma.

\proof {\sc of Theorem \ref{lbez}}
By the previous Lemma, part one, 
\[ v_{(x,x)}(Y \# Z) \geq v_x(Y) v_x(Z). \] 
Further, by Remark \ref{luraum}, 
\[ v_{(x,x)}(Y \# Z) \leq v_{(x,x)} (\Pe(\Delta). (Y \# Z)) = 
   v_{(x,x)} (\delta_*(Y . Z)), \]
which by part 2 of the previous Lemma equals $v_x(Y.Z)$.

\satz{Definition}
Let $\CY$ be an effective cycle in $\Pe^M_\Z$, and
$Y$ its base extenseion to $\spec \; \Q$. For a real number $H$ 
the weighted order of vanishing of $\CY$ at a point $x$ in $\Pe^M_k$ is 
defined as $v_x(Y)/t_H(\CY)$.
\end{Satz}

\satz{Lemma} \label{wov}
For every effective cycle $\CY$, and every closed point $x \in \Pe^M$,
there is an irreducible component $\bar{\CY}$ of $\CY$ such that
\[ \frac{v_x(\bar{Y})}{t_H(\bar{\CY})} \geq \frac{v_x(Y)}{t_H(\CY)}. \]
\end{Satz}

\proof
Follows from the fact that both $v_x$ and $t_H$ are linear functions
on $Z(\Pe^t)$, and elementary arithmetic.

\satz{Proposition} \label{mult}
Let $X \subset \Pe^M$ be an irreducible subvariety of dimension $t$, and 
$\alpha$ a closed point in $X$. Further, $Y$ a subvariety of codimension 
$p$ in $X$ containing $\alpha$, and
$f_i \in \Gamma(\Pe^M,O(D_i)), i=1, \ldots t-p$ global sections such that
for every $i = 0, \ldots, t-p$ there is an effective cycle $Z_i$ of pure
codimension $i+p$ such that $Z_0 = Y$, the intersection of $Z_i$ with
$\di f_{i+1}$ is proper, and $Z_{i+1} + X_i  = \di f_{i+1} . Z_i$, where
$X_i$ is an effective cycle whose support does not contain $\alpha$.
Further, assume that
for every $i=1, \ldots, t-p$ there is a number $S_i \in \N$ 
such that $\partial^I f_i$ is zero on 
$\alpha$ for every $i=1, \ldots, t-p$, $I$ with $|I| \leq S_i$, and
$\partial^I$ a derivation of the functions field of $X$ as above. Then, 
\[ v_\alpha(Z_{t-p)} \geq  S_1 \cdots S_{t-p}. \]
\end{Satz}

\proof
By fact \ref{dim0}.1, $v_\alpha(Y) \geq 1$, and
by Proposition \ref{vfunc}, the vanishing order of $f_i$ at $\alpha$ is
at least $S_i$. Hence, by the local B\'ezout Theorem,
\[ v_\alpha(Z_{i+1}) = v_\alpha(Z_{i+1} + X_i) = v_\alpha(\di f_{i+1} . Z_i) =
   v_\alpha((X. \di f_{i+1}). Z_i) \geq \]
\[ v_\alpha(X) v_\alpha(\di f_{i+1}) v_\alpha( Z_i) \geq 1 S_{i+1} v_\alpha(Z_i), \]
and the Proposition follows by complete induction.

%\satz{Definition} \label{sucint}
%Let $f_1, \ldots, f_n$ be global sections on $\Pe^t$. They are said to
%sucessively intersect properly, if there is a chain
%\[ \CY_0 = \Pe^t \supset \CY_1 \supset \cdots \supset \CY_n \]
%such that for every $k=1, \ldots, n$ the intersection $\di f_k . \CY_{k-1}$
%is proper, and $\CY_{k}$ is the irreducible component in 
%$\di f_k . \CY_{k-1}$ that has minimal weighted derivated algebraic distance. 
%\end{Satz}

%\satz{Proposition}
%Let $\CX$ be an irreducible projective variety of relative dimension $t$, 
%$\bar{\CL}$ a positive metrized very ample line bundle, 
%$\CX \hookrightarrow \Pe^M$ the embedding defined by $\CL$, and
%$\theta \in X(\C)$ a generic point. Further
%$\Pe^t \subset \Pe^M$ a projective subpace such that the projection $\pi$
%of $\Pe^M$ to $\Pe^t$ restricted to $\CX$ is   

%\begin{enumerate}

%\item
%If $\CY$ is an effective cycle of pure codimension $p$ in $\CX$,
%that can be written as the proper intersection of a cycle $\bar{\CY}$
%on $\Pe^M$, and $\CX$, then
%\[ \deg_{O(1)} \pi(Y) = \deg_L X \deg_L \CY, \]
%\[ h_{\overline{O(1)}} (\pi(\CY)) \leq \deg_L X h(\CY) +
%   \deg Y h(\CX) + c \deg X \deg Y, \]
%and for $S \leq \deg Y/3$,
%\[ D^S(\pi(\theta),\pi(Y)) \leq D(\theta,Y) + 
%   O(\deg (X.Y) + h(\CX. \CY) \log (\deg (X.Y) + h(\CX. \CY))). \]

%\item
%If additionally $\CY = \di f$ with $f$ a global section of $\CL^{\otimes D}$,
%then $X f = \pi X \pi f$.
%\end{enumerate}
%\end{Satz}

\subsection{Weighted derivative algebraic distance}

In analogy to the weighted algebraic distance defined in \cite{App3},
define the weighted derivated algebraic distance

\satz{Definition}
Let $\CX$ be an effective cycle in $\Pe^t$. The $a$-size of $\CX$ is 
defined to be the number
\[ t_a(\CX) := a \deg X + h(\CX). \]
For $S \in \N$ define the weighted derivated algebraic distance of $\CX$ to
$\theta$ as
\[ \varphi_a^S(\theta,\CX) := \frac{D^{3S}(\theta,X)}{t_a(\CX)}. \]
\end{Satz}

\satz{Lemma} \label{gaal}
Let $\CX$ be an effective cycle in $\Pe^t$, and $S$ a natural number. Then,
there is an irreducible component $\CY$ of $\CX$ and an $S_Y \in \N$ with
$S_Y/t_a(\CY) \geq S/t_a(\CX)$ such that
\[ 2 \varphi_a^{S_Y}(\theta,\CY) \leq \varphi_a^S(\theta,\CX) +
   O\left( \frac{\log \deg X}a \right). \]
We call $\CY$ the irreducible component with minimal derivated
algebraic distance relative to $S$.
\end{Satz}

\proof
Let $\CY, \CZ$ be effective cycles of codimension $p$ in $\Pe^t$, and 
$S \in \N$. By \cite{App4}, Theorem 5.1, there are subspaces 
$\Pe(F), \Pe(F')$ of codimension $t-p$ such that
with $y_1, \ldots y_{\deg Y}$ the points in the intersection of 
$\Pe(F)$ with $Y$ counted with multiplicity, and likewise
$z_1, \ldots, z_{\deg Z}$ for $\Pe(F')$ and $Z$ for all natural numbers
$S_1 \leq \deg Y/3, S_2 \leq \deg Z/3$,
\[ D^{S_1}(Y,\theta) \leq \sum_{i=S_1+1}^{\deg Y} \log |y_i, \theta| +
   O((S_1+\deg Y) \log (S_1 \deg Y)), \]
\[ D^{S_2}(Z,\theta) \leq \sum_{i=s_2+1}^{\deg Z} \log |z_i, \theta| +
   O((S_2+\deg Z) \log (S_2 \deg Z)). \]
For given $S$ choose $S_1,S_2 \in \N$ such that $S_1+S_2=S$, and
\[ \log |y_i,\theta| \leq \log |z_j,\theta| \quad \mbox{for} \quad
   i \leq S_1, j \geq S_2, \]
\begin{equation} \label{gaal1}
   \log |z_j,\theta| \leq \log |y_i,\theta| \quad \mbox{for} \quad
   j \leq S_2, i \geq S_1. 
\end{equation}
Then,
\[ \sum_{i=S_1+1}^{\deg Y} \log |y_i,\theta| +
   \sum_{j=S_2+1}^{\deg Z} \log |z_j,\theta| \leq
   \inf_{\Pe(F)} \sum_{z \in supp(\Pe(F).(X+Z)} n_z \log |z,\theta| \leq \]
\[ \frac12 D^{3S}(\theta,X+Y) + O(\deg (X+Y) \log (\deg (X+Y))), \]
again by \cite{App4}, Theorem 5.1.
Consequently,
\[ \frac{\sum_{i=S_1+1}^{\deg Y} \log |y_i,\theta| +
   \sum_{j=S_2+1}^{\deg Z} \log |z_j,\theta|}{t_a(\CX+\CY)} \leq
   \frac12 \frac{D^{3S}(\theta, Y+Z)}{t_a(\CY+\CZ)} + 
   O\left(\frac{\log (\deg (Y+Z))}a\right). \]
Let $r \in \R$ be such that 
$\frac{S_1 +r}{t_a(\CX)} = \frac{S}{t_a(\CX + \CY)}$,
and $s = \mbox{sign} r [|r|]$. Then, 
$|s| \leq \mbox{min}(S-S_1,S-S_2)$, and by (\ref{gaal1}),
\[ \frac{\sum_{i=S_1+1}^{\deg Y} \log |y_i,\theta| +
   \sum_{j=S_2+1}^{\deg Z} \log |z_j,\theta|}{t_a(\CY+\CZ)} \geq \]
\[ \frac{\sum_{i=S_1+s+1}^{\deg Y} \log |y_i,\theta| + 
   \sum_{j=S_2-s+1}^{\deg Z} \log |z_j,\theta|}{t_a(\CY)+t_a(\CZ)}. \]
By elementary arithmetic, this is greater or equal
\[ \mbox{min} 
   \left(\frac{\sum_{i=S_1+s+1}^{\deg Y} \log |y_i,\theta|}{t_a(\CY)},
     \frac{\sum_{j=S_2-s+1}^{\deg Z} \log |z_j,\theta|}{t_a(\CZ)} \right). \]
Further 
\[ \frac{S_1+r}{t_a(\CX)} = \frac{S_2 -r}{t_a(\CY)} = \frac{S}{t_a(\CX+\CY)}.\]
By complete induction it follows, that for any effective cycle
$\CX$ with decomposition into irreducible parts
\[ \CX = \sum_{k=1}^n \CX_1, \]
we have numbers $S_k$ with $S_k /t_a(\CX_k) = S/t_a(\CX)+ \epsilon$, and
\[ \mbox{min}_{k=1, \ldots n} \left(
   \frac{\sum_{i=S_k+1}^{\deg X_k} |\log x_{ik},\theta|}{t_a(\CX_k)}\right)\leq
   \frac{D^{3S}(\theta,X)}{2t_a(\CX)} + O((\log \deg X)/a). \]
The Lemma follows by once more using \cite{App4}, Theorem 5.1.

\satz{Lemma} \label{gaaleins}
Let $\CY \in Z_{eff}(\Pe^M)$ an effective cycle, and 
$\theta \in \Pe^M(\C)$ a point not contained in the support of $\CY$.
Then, for any $m,n,S \in \N$.
\[ D^{nS} (\theta,m n Y) \leq m D^{nS}(\theta,Y). \]
\end{Satz}

\proof
Since
\[ \exp(D(\theta,mnY)) = \left(\exp(D(\theta,X)) \right)^{mn}, \]
this follows by elmentary differentiation techniques.

\section{Projection to a projective sub space} \label{psp}

\satz{Proposition}
Let $X \subset \Pe^M_\C$ be a subvariety of dimension $t$, further
$\theta \in X(\C)$, and $Y$ an effective cycle in $X$ whose support does
not contain $\theta$.
Let $\varphi: \A^M(\C) \to \Pe^M(\C)$ be an affine chart centered at $\theta$
such that $\varphi(\A^t \times \{0\} = \Pe(T_\theta X)$ the tangent space
of $X$ at $\theta$. Denote by $I$ a multi index, and by
$\partial^I$ the corresponding differential. Further let $N_t$ be 
the set of multi indizes $I =(1_1, \ldots, i_{2M})$ with
$i_{2t+1}= \cdots = i_{2M} = 0$. Then, for $S \leq \deg Y/3$,
\[ \sup_{I \in N_t, |I| \leq S} \log 
   \left| (\partial^I (\varphi^* \exp D(Y,\theta)) \right| \leq
   D^S(\theta,Y). \]
\[ D^S(\theta,Y) \leq \sup_{I \in N_t, |I| \leq S} \log 
   \left| (\partial^I (\varphi^* \exp D(Y,\theta)) \right| +
   O(\deg Y \log \deg Y). \]
\end{Satz}

%Since $\theta$ is a generic point of $X$ the orthogonal projection
%from $X(\C)$ to the tangent space of $X$ at $\theta$ is an affine chart, 
%and the claim follows from \cite{App4}, Lemma 3.1.

\proof
The first claim is trivial. For the second claim, let $U_\theta$ be a 
neighbourhood of $\theta$ in $X$ such that the orthogonal projection
$\pi$ of $U_\theta$ to $T_\theta X$ is bijective, and for every 
$x \in U_\theta$, the inequality $|x,\theta| \leq 2 |\pi x, \theta|$ holds.
With $p$ the codimension of $Y$ in $X$, by \cite{App4}, Theorem 1.4, 
there is a subspace
$\Pe(F) \subset \Pe^M$ of codimension $t-p$ such that 
$\Pe(F)$ contains $\theta$, intersects $Y$ properly, and with
$\Pe(F). Y = \sum_{i=1}^{\deg Y} y_i$, numbered in such a way that
$|y_1,\theta| \leq \cdots \leq |y_{\deg Y},\theta|$ the derivated algebraic
distance of $\theta$ to $Y$ may be estimated as
\[ D^S(\theta,Y) \leq \sup_{|I| \leq S} \log |\partial^I 
   \prod_{i=1}^{\deg Y} |y_i,\theta|| + O(S \log \deg Y), \]
\[ \log |\partial^I \prod_{i=1}^{\deg Y} |y_i,\theta|| \leq 
   D^S(\theta,Y) + O(\deg Y). \]
Let $r$ be the radius of $U_\theta$, and $k \leq \deg Y$ a number such that
$|y_k,\theta| \leq r \leq |y_{k+1},\theta|$. Then, with 
$c_i(z) = |\pi y_i,\theta|/|y_i,\theta|$, 
\[ \log |(\partial^I \varphi^*c_i(z))(0)| \leq c |I|, \quad
   \log |(\partial^I 1/c_i(z))(0)| \leq c |I|, \]
with $c$ a fixed constant. Hence, for every $I$, with $|I| \leq S$,
\[ \log \sup_{|I| \leq S}|\partial^I\prod_{i=1}^{\deg Y} |y_i,\theta|| \leq 
   \log \sup_{|I| \leq S}|\partial^I \prod_{i=1}^{\deg Y} |\pi y_i, \theta||
   + c S \leq \]
\[ \log \sum_{I \in N_t, |I| \leq S} \partial^I \prod_{i=1}
   |\pi y_i,\theta|| + c S \leq
   \log \sup_{I \in N_t, |I| \leq S} |\partial^I (\varphi^* \exp D(Y,\theta))|,
   \]
proving the second claim.

\satz{Lemma} \label{reda}
There are positive constants $\bar{c}, \tilde{c}$ only depending on $M$ and 
$t$, and a subspace 
$\Pe^{M-t-1} \subset \Pe^M$ defined over $\Z$ that does not meet $\CX$, and
fulfills
\[ h(\Pe^{M-t-1}) \leq \tilde{c} \log \deg X \quad \mbox{and} \quad
   \log |\Pe^{M-t-1},X| \geq - \bar{c} -\log \deg X. \]
For $\Pe^t$ the orthogonal complement of $\Pe^{M-t-1}$ in $\Pe^M$, the 
restriction of the map
\[ \pi: \Pe^M \setminus \Pe^{M-t-1}, \quad [v+w] \mapsto [v], 
   [v] \in \Pe^t, [w] \in \Pe^{M-t-1} \]
to $\CX$ is a flat, finite proper map $\pi_X: \CX \to \Pe^t$, and
\[ h(\Pe^t) \leq c \log \deg X, \]
with $c$ a constant only depending on $M$ and $t$.
\end{Satz}

\proof
By \cite{App4}, Corollary 5.4, there
is a subspace $\Pe(W) \subset \Pe^M_\C$ with
\[ \log |\Pe(W),X| \geq - c_1 - \log \deg X, \]
with some positive constant $c_1$ only depending on $M$, and $t$.
For $V \subset \C^{M+1}$ a subspace, denote by $S(V)$ the set of vectors
of length $1$ in $V$, and by $pr_V^\bot$ the orthogonal projection to the
orthogonal complement of $V$. On the Grassmannian $G_{M,t,}$ we have
\[ |V,W| = \sup_{v \in S(V)} |pr_W^\bot|, \]
and for $V$ a primitive submodule of $\Z^{t+1}$, 
\[ h(V) = \log \mbox{vol} (V)+ \sigma_p. \]
Let $W$ be the space from above, $q=M-t = \dim W$, and 
$a=2 e^{c_1} q (t+1)\deg X$. One can recursively 
find vectors
\[ v_1, \ldots, v_q \in \C^{M+1}, \]
such that with $V_i = \la v_i, \ldots, v_i \ra$,
\[ v_i \in pr^\bot_{V_{i-1}}(\Z^{M+1}), \quad 
   |v_i| \leq (M+1) 2^{(M+1)/q} a^{t+1}, \quad
   |pr_W^\bot(v_i)| \leq \sqrt{t+1} \frac1a. \]
Indeed, assume that $w_1, \ldots, w_{t+1}$ is an orthonormal basis
of $W^\bot$, and $v_1, \ldots, v_i$ have been found. Since,
$\log \mbox{vol} V_i \geq 1$, then $\log \mbox{vol}(\Z^{M+1}/V_i)) \leq 1$.
Let $Q$ be the Cuboid in $\R^{t+1}$ that has lengths
$2^{(M+1)/(t+1)} a^{q/(t+1)}$ parallel to $W$, and lengths
$1/a$ paralell to $W^\bot$. Then,
\[ \mbox{vol} (Q) =
   2^{M+1} a^{t+1} (1/a)^{t+1} =
   2^{M+1} \geq 2^{M+1} \mbox{vol}(\Z^{M+1}/V_i). \]
By the Theorem of Minkovksi, $Q$ thus contains a non zero vector $v_{i+1}$, 
and $v_{i+1}$ fulfills
\[ |v_{i+1}|^2 \leq q (2^{(M+1)/(t+1)} a^{q/(t+1)})^2 + (t+1) (1/a)^2 \leq
   (M+1) 2^{(M+1)/(t+1)} a^{2q/(t+1)}, \]
and 
\[ |pr_W^\bot(v_{i+1})|^2 \leq (t+1) \left(\frac1a\right)^2. \]
Since $v_1, \ldots, v_q$ is an orthonormal basis of $V = V_q$,
for any $v \in S(V)$ we have $v= \sum_{i=1}^q a_i v_i$ with
$\sum_{i=1}^q |a_i|^2 =1$, hence
\[ |pr_W^\bot(v)| \leq \sum_{i=1}^q a_i |pr_W^\bot(v_i)| \leq
   q \sqrt{t+1} \; \frac1a = \frac{q (t+1)}{2 e^{c_1} q(t+1) \deg X} =
   \frac{e^{-c_1}}{2 \deg X}. \]
Hence, 
\[ \log |V,W| = \log \sup_{v \in S(V)} |pr_W^\bot(v)| \leq - c_1 - \log 2
   - \log \deg X. \]
Since $\log |W,X \geq - c_1 - \log \deg X$, we get
$\log |V,X| \geq - c_1 - \log \deg X- \log 2 = - \bar{c} - \log \deg X$ with
a suitable $\bar{c}$.

Finally, since $|v_i| \leq (M+1) 2^{(M+1)/q} a^{t+1}$ for $i=1, \ldots q$,
we get
\[ h(\Pe(V)) = \log \vol(V) + \sigma_q  
   \sum_{i=1}^q  \log |v_i| + \sigma_q \leq \]
\[ \log \left(q (M+1) 2^{(M+1)/(t+1))} 
   (2 e^{c_1} q (t+1) \deg X)^{q/(t+1)} \right) + 
   \sigma_q \leq \tilde{c} \log \deg X, \]
with a suitable $\tilde{c}>0$.
If $M = \Z^{M+1} \cap V$, and $M^\bot = \Z^{t+1} \cap V^\bot$, by
\cite{Be}, Proposition 1.(ii),
\[ \mbox{vol} M^\bot \mbox{vol} M \leq (\deg X)^{\tilde{c}}/\exp(\sigma_q). \]
Hence, with $\Pe^t = \Pe(M^\bot)$,
\[ h(\Pe^t) = \log \mbox{vol} M^\bot + \sigma_t \leq
   \tilde{c} \log \deg X + \sigma_t - \sigma_q \leq c \log \deg X. \]

\satz{Proposition} \label{runter}
Let $\CY \in Z_{eff}^p(\CX)$ be a cycle, $\theta \in X(\C)$ 
a generic point, and $\Pe^t,\Pe^{M-t-1},\pi, \pi_X$ as in 
Lemma \ref{reda}

\begin{enumerate}

\item
If the set of complex valued points $Y_i$ of an irreducible component $\CY_i$ 
of $\CY$ has sufficiently small distance to $\theta$, then
$\dim \pi(\CY_i) = \dim \CY_i$.

\item
For $x,y \in X(\C)$ in a sufficiently small neighbourhood of $\theta$,
\[ |x,y| \leq c |\pi_X x, \pi_X y|, \]
where $c$ is constant depending on $\theta$.
and for $x,y \in \Pe^M \setminus \Pe(F^\bot)$,
\[ \log |\pi x, \pi y| \leq |x,y| -
   \log |x,\Pe^{M-t-1}| - \log |y,\Pe^{M-t-1}|. \]

\item
If $Y$ is irreducible, $\dim \pi Y = \dim Y$, then
\[ \deg \pi_X(Y) = \deg Y, \quad h(\pi_X(\CY)) \leq  h(\CY). \]

\item
If $Y$ is irreducible, $\dim \pi Y = \dim Y$,
$\theta$ is not contained in the support of $Y$, and 
$S \leq \deg Y/3$, then
\[ 2 (D^{\Pe^t})^S(\pi \theta, \pi_X(Y)) \leq D^{3S}(\theta, Y) + 
   O(\deg Y \log \deg Y), \]

%\item
%For every $f \in \Gamma(\CX,\CL^D)$ such that $f(\theta)$ is sufficiently
%small, there is an $F \in \Gamma(\Pe^t,O(D))$ such that
%$\di(\pi_X^*F)- \di f$ is an effective cycle, 
%\[ \log |F| \leq \log |f| + D h(\CX), \]
%and with
%$\partial^I$ the differential operators from Theorem \ref{Abl}, and
%$\tilde{\partial}^I$ the canonical affine differentail operator on
%$\Pe^t$,
%\[ \sup_{|I| \leq S} \log |(\tilde{\partial}^I F)(\pi \theta)| \leq
%   \sup_{|I| \leq S} \log |(\partial^I f) (\theta)| + O(S). \]
\end{enumerate}
\end{Satz}

\proof
1. Let $T_\theta X$ be the tangent space of $\CX$ at $\theta$ which may be
identified with the projective space $\Pe(V_\theta)$ corresponding to a 
subspace $V_\theta$ of $\C^{M+1}$. Since $\theta$ is a generic point of $\CX$, 
the restriction of $\pi$ to $\Pe(V_\theta)$ is bijective, hence comes from
a bijective linear map
\[ \varphi: V_\theta \to \C^{t+1}. \]
Because the metrics on $\Pe(V_\theta)$ and $\Pe^t$ just correspond to different
inner products on $V_\theta$ and $\C^t$, there is positive constant $c$ such that
\[ \frac1c |\pi x, \pi y| \leq |x,y| \leq c |\pi x, \pi y| \]
for every $x,y \in \Pe(V_\theta)$.
Further, for a sufficiently small neighbourhood $U_\theta$ of $\theta$ the
orthogonal projection $pr$ from $U_\theta$ to $T_\theta X=\Pe(V_\theta)$ is 
bijective, and
\[ |pr x, pr y| \leq |x,y| \leq 2 |pr x, pr y| \]
for every $x,y$ in $U_\theta$ implying the first claim.

Let $u,v \in \C^{t+1}$ be vectors representing $\pi x$ and $\pi y$.
There are vectors $w_1, w_2 \in \C^{M-t}$ such that
$\bar{u} = u + w_1, \bar{v} = v + w_2$ represent the points
$x,y$. We may assume that $|\bar{u}| = |\bar{v}| =1$. Then, in the 
Fubini-Study metric, since $u,w \in \C^{t+1}$, and 
$w_1,w_2 \in \C^{M-t} = (\C^{t-1})^\bot$,
\[ |x,\Pe^{M-t-1}|^2 \leq |x,[w_1]|^2 = \sin^2 (u,w_1) = |u|^2, \quad
   |y,\Pe^{M-t-1}| \leq |y,[w_2]|^2 = |v|^2. \]
Without loss of generality, we may assume $\la u |v\ra \leq 0$, and
$|u| \leq |v|$, hence $|w_2| \leq |w_1|$. 
If $|w_2| = 0$, then $x=\pi x$, $y = \pi y$, and there is nothing to prove.
If $|w_2|>0$, set $\lambda = |w_2|/|w_1| \leq 1$, and define the point 
$\tilde{y} \in \Pe^M$ by $y = [v + \lambda w_1]$. 

Then,
\[ |x,y|^2 = 1- (\la u | v \ra + \la w_1 | w_2 \ra)^2 \geq
   1 - (\la u | v \ra + |w_1| |w_2|)^2 = \]
\[ 1 - (\la u | v \ra + \lambda \la w_1 | w_1 \ra)^2 =
   |x,\tilde{y}|^2. \]
Further,
\[ |u|^2 |v|^2 |\pi x| \; |\pi y|^2 =
   |u|^2 |v|^2 (1- \la u | v \ra^2) \leq \]
\[ |u|^2+|v|^2-|u|^2 |v|^2 - 2 \la u |v \ra |w_1| |w_2| - \la u|v\ra^2 = \]
\[ 1- (1- |u|^2)(1-|v|^2) - \la u |v\ra^2 - 2 \la u |v\ra |w_1| |w_2| =
   1- \la u | v \ra^2 - 2 \la u |v\ra \lambda |w_1|^2 - \lambda^2 |w_1|^4 = \]
\[ 1 - \la u|v\ra + \la w_1|\lambda w_2 \ra^2 = |x, \tilde y|, \]
which, together with the above, implies
\[ |x,\Pe^{M-t-1}|^2 |y,\Pe^{M-t-1}|^2 |\pi x| \pi y|^2 \leq
   |x,y|. \]

\vspace{2mm}

2. Since $\theta$ is a generic point, the base extension $\pi_\C$ to $X_\C$ 
is injective in some neighbourhood of $\theta$. This immediately implies the 
claim.

\vspace{2mm}

3. The first claim is obvious. The second claim is \cite{BGS}, (3.3.7).

\vspace{2mm}

4. Let $p$ be the dimension of $Y$.
Since $|\Pe^{M-t-1},X| \geq - \bar{c} - \log \deg X$, by \cite{App4},
Propositions 5.4, and Corollary 5.5, there is a space $\Pe(F) \subset \Pe^M$
of codimension $t-p$ that contains $\Pe^{M-t-1}$ as well as $\theta$, hence
intersects $Y$ properly, such that
\[ (D^{\Pe(F)})^S(\theta,Y. \Pe(F)) \leq D^S(\theta,Y) + 
   O(\deg Y \log \deg Y), \]
hence, if $\Pe(F).Y = \sum_{i=1}^{deg Y} y_i$ where the $y_i$ are ordered in 
such a way that $|y_1,\theta| \leq \cdots \leq |y_{\deg Y},\theta|$, \cite{App4},
Proposition 4.7 implies
\[ 2 \sum_{i=S+1}^{\deg Y} \log |y_i,\theta| \leq D^{3S} (\theta,Y) +
   O(\deg Y \log \deg Y). \]
Let $\sigma \in \Sigma_{\deg Y}$ be a permutation such that
$|\pi y_{\sigma 1},\pi \theta| \leq\cdots\leq |\pi y_{\sigma \deg Y},\theta|$. 
By part Proposition \ref{hoch}, $|\pi y_i, \pi \theta| \leq |y_i, \theta| + c \log \deg Y$. Hence,
\[ 2 \sum_{i=S+1}^{\deg Y} \log |\pi y_{\sigma i},\pi \theta| \leq
   2 \sum_{i=S+1}^{\deg Y} \log |y_{\sigma i},\theta| + c(\deg Y-S)
   \log \deg Y \leq \]
\[ 2 \sum_{i=S+1}^{\deg Y} \log |y_i,\theta| + c (\deg Y-S) \log \deg Y \leq
   D^{3S} (\theta,Y) + O(\deg Y \log \deg Y). \]
Further, since $\Pe(F) \cap \Pe^t$ is a subspace of dimension $p$ in $\Pe^t$
containing $\pi \theta$ and intersecting $\pi Y$ properly,
\cite{App4}, Proposition 5.1 implies
\[ (D^{\Pe^t})^S(\pi \theta,\pi Y) \leq \sum_{i=S+1}^{\deg Y} 
   \log|\pi y_{\sigma i},\pi \theta| + O(\deg Y \log \deg Y), \]
hence
\[ 2 (D^{\Pe^t})^S(\pi \theta,\pi Y) \leq D^{3S} (\theta,Y) + 
   O(\deg Y \log \deg Y), \]
as was to be proved.

\vspace{2mm}

\satz{Lemma} 
Let $\CY \in Z_{eff}^p(\Pe^M)$ be an effective cycle that intersects
$\CX$ properly,
and $\theta \in X(\C)$ a point not contained in the support of $Y$. Then,

\begin{enumerate}
\item
\[ \deg (X . Y) = \deg X \deg Y, \]
\[ h(\CY) \leq \deg h(\CY) + \deg Y h(\CX) + c \deg X \deg Y. \]

\item
For any $S \leq \deg Y$,
\[ 2 D^S(\theta, Y . X) \leq D^(3S)(\theta,Y) + 
   O(\deg X \deg Y \log(\deg X \deg Y)) . \]
\end{enumerate}
\end{Satz}

\proof
1. is just the algebraic and arithmetic B\'ezout Theorem. 
Since $\theta \in X(\C)$,
2. is Theorem \ref{DMBTcor}.2 applied to the varieties $X,Y$.

\satz{Proposition} \label{hoch}
In the situation of Lemma \ref{reda}, let $\CY \in Z_{eff}^p(\Pe^t)$.
Then, $\CX$ intersects $\pi^* (\CY)$ properly, and
$\CY^* := \pi^*_{\CX}(\CY) = \pi^*(\CY) . \CX$. Further,
\begin{enumerate}

\item
\[ \deg Y^* = \deg X \deg Y, \]
\[ h(\CY^*) \leq \deg X (h(\CY) + \tilde{c} \deg Y \log \deg X) +
                   \deg Y h(\CX) + c \deg X \deg Y, \]
and for every irreducible component $\bar{\CY}^*$ of $\CY^*$ sufficiently
close to $\theta$,
\[ \deg \bar{Y}^* \geq \deg Y, \quad h(\bar{\CY}^* \geq h(\CY). \]

\item
If further $\theta \in \Pe^t(\C)$ is not contained in the support
of $Y$, and $\bar{\theta} \in \CX(\C)$ is a point with 
$\pi_X \bar{\theta} = \theta$, then for $S \leq \deg Y$,
\[ D^S(\bar{\theta}, Y^*) \leq \frac14 D^{9S}(\theta, Y) +
   \deg X h(\CY^*) + \deg Y^* h(\CX) + d \deg X \deg Y^*. \]

\item
If $f \in \Gamma(\Pe^t,O(D))$, let $f^* \pi^* f$.
Then,
\[ \log |f^*|_{L^2(\Pe^M)} = |f|_{L^2(\Pe^t)} + c D, \]
\[ |\di f^*, \theta| \leq c|\di f, \pi \theta| \leq c|\di f^*, \theta| +
   c c_2 \deg X. \]
\[ \sup_{|I| \leq S} \log |(\partial^I f^*)(\theta)| \leq
   \sup_{|I| \leq S} \log |(\partial^I f) (\bar{\theta})|. \]

\end{enumerate}
\end{Satz}

\proof

1. Since $\deg \pi^* Y = \deg Y$, the first claim follows from
$\pi^*_{\CX}(\CY) = \pi^*(\CY) . \CX$ and the previous Lemma.

Let $x_1, \ldots, x_{M-t} \in \Gamma(\Pe^M,O(1))$ such that
$\Pe^t = \di x_1 . \ldots . \di x_{M-t}$. Then, by Lemma \ref{bashoeh},
\[ \sum_{i=1}^{M-t} \int_{\di x_1 . \ldots . \di x_{i-1}} 
   \log |x_i| \mu^{M-i} = h(\Pe^t) - h(\Pe^M), \]
and $\CY = \pi^* (\CY) . \di x_1 . \ldots . \di x_{M-t}$.
Hence, there are numbers $a_1, \ldots, a_{M-t} \in \R$ such that
$\sum_{i=1}^{M-t} a_i = h(\Pe^t)- h(\Pe^M)$, and 
$\log |x_i| - a_i$ is a normalized Green form for 
$\di x_1 . \ldots . \di x_i $ in $\di x_1 . \ldots . \di x_{i-1}$. The
equality $\CY = \pi^*(\CY) . \di x_1, \ldots . \di x_{M-t}$ together with
Lemma \ref{bashoeh} and \cite{BGS}, Proposition 5.1 implies
\begin{eqnarray*}
h(\CY)- h(\pi^*(\CY)) &=&
   \sum_{i=1}^{M-t} \int_{\pi^*(Y) . \di x_1 . \ldots . \di x_{i-1}} 
   \log |x_i| \mu^{m-p-i} \\
  &=&\sum_{i=1}^{M-t} \int_{\pi^*(Y) . \di x_1 . \ldots . \di x_{i-1}}
   (\log |x_i| -a_i) \mu^{M-p-j}  + \deg Y \sum_{i=1}^{M-t} a_i \\
&=& -c \deg Y
   + \deg Y (h(\Pe^t) - h(\Pe^M), 
\end{eqnarray*}
with $c$ a positive constant depending only on $t,M$, and $p$.
Thus, 
\[ h(\pi^*(\CY)) = h(\CY) + c \deg Y - \deg Y (h(\Pe^M) - h(\Pe^t))  \leq \]
\[ h(\CY) + c \deg Y + c_1 \deg Y \log \deg X. \] 
Since $\pi^*_{\CX}(\CY) = \pi^*(\CY) . \CX$, the previous Lemma implies
\[ h(\pi^*_{\CX}(\CY)) \leq \deg X h(\pi^*(\CY)) + \deg Y h(\CX) + 
   c_2 \deg X \deg Y \leq \]
\[ \deg X (h(\CY) + c_1 \deg Y \log \deg X) + \deg Y h(\CX) + 
   c_3 \deg X \deg Y, \]
proving the second claim.

If $\bar{\CY}^*$ is an irreducible component of $\CY^*$ sufficiently close
to $\theta$, then because of the irreduciblity,
$(\pi_X)_* \bar{\CY}^* = \CY$, hence by Proposition \ref{runter}.2,
$\deg \bar{Y}^* \leq \deg Y, h(\bar{\CY}^*) \leq h(\CY)$.

\vspace{2mm}

2. Let $U_{\theta}$ be a sufficiently small neighbourhood of $\theta$ in
$X(\C)$.

By \cite{App4}, Theorem 1.4, there is a subspace $\Pe(F) \subset \Pe^t$ of 
dimension  $p$ such that with $\Pe(F) . Y = \sum_{i=1}^{\deg Y} y_i$,
ordered such that $|y_1,\theta| \leq \cdots \leq |y_{\deg Y},\theta|$,
\[ 2 \sum_{S+1}^{\deg Y} 
   \log |y_i,\theta| \leq D^{3S}(\theta,Y) + O((S+\deg Y) \log \deg Y), \]
for every $S \leq \deg Y/3$.
Denote by $l \leq \deg Y$ the number such that
$y_i \in \pi_X U_\theta$ for $i \leq l$, and $y_i \notin \pi U_\theta$ 
for $i > l$. Then, $\log |y_i,\theta| \geq -c_2$ for every $i>l$
with $c_2>0$ independent of $Y$.
Further, let $\Pe(F^*) \subset \Pe^M$ be the projective 
subspace of codimension 
$p$ that Contains $\Pe(F)$ as well as $\Pe^{M-t-1}$. Then the restriction
of $\pi_X$ to $U_\theta$ maps
$\Pe(F^*) \cap  \mbox{supp}(\pi^*(Y))$ injectively to 
$\Pe(F \cap \mbox{supp} Y)$, and
since $|\bar{\theta}, \Pe^{M-t}| \geq c \deg X$, for every $y^*$ in
$\Pe(F^*) \cap \pi^*(Y)$, we have
$\log |y^*,\bar{\theta}| \leq \log |\pi(y^*),\theta| + c_1 \log \deg X$, and
consequently if $\Pe(F^*) . \pi^*(Y) = \sum_{i=1}^{\deg Y} y^*_i$ ordered
in the usual way,
\begin{eqnarray*}
 D^S(\pi^*Y,\theta) &\leq&
   \sum_{i=S+1}^{\deg Y} \log |y_i^*,\theta| + O( S \log \deg Y) \\
&\leq& \sum_{i=S+1}^l \log |y_i^*,\theta| + O(S \log \deg Y) \\
&\leq& \sum_{i=S+1}^l \log |y_i,\theta| + 
   c_1 \deg Y \log \deg X +O(S \log \deg Y) \\
&\leq& \sum_{i=S+1}^{\deg Y} \log|y_i,\theta| + (c_2+c_1) \deg Y \log \deg X +
   O(S \log \deg Y). 
\end{eqnarray*}
Hence,
\[ D^S(\pi^*(Y), \bar{\theta}) \leq
   \frac12 D^{3S} (Y,\theta) + (c_2+c_1) \deg Y \log \deg X +O(S \log \deg Y).
\]

\vspace{2mm}

3.
The first claim follows by integration over the fibres of $\pi$, and the
second claim from Proposition \ref{hoch}.1.

With $\varphi: \A^t \to \Pe^t$ the canonical affine chart centered at
$\theta$, and $\psi$ the local inverse of $\pi_X$ at $\bar{\theta}$ with
image in $U_:\theta$, the map $\psi \circ \varphi$ is an affine chart
of $\CX$ around $\theta$. Thus, for an $f \in \Gamma(\Pe^t,O(D))$,
\[ (\psi \circ \varphi)^* \circ \pi^*f = \varphi^*f, \]
from which the claim about derivatives follows.

The inequality
$|\di f^*, \theta| \leq |\di f, \pi \theta| \leq |\di f^*, \theta| +c \deg X$
follows from part 1.

\section{Proof of second criterion}

This section establishes a proof of
Theorem \ref{algind1}. For a given $a > 1$, if $H_k \leq a D_k$,
one can replace $H_k$ by $\bar{H}_k = a D_k$. Then, since
$\bar{H}_k + D_k \leq (a+1) D_k \leq (a+1) (D_k + H_k)$, still
\[ \limsup_{k \to \infty} \frac{S_k^s V_k}{D_k^s (D_k + \bar{H}_k)}= \infty, \]
hence we may from now on assume that $H_k \geq a D_k$. For similar reasons, one
may assume $S_k \leq 3 D_k$ for all $k$.
Similarly, by replacing the series $(D_k,H_k,S_k,V_k)$ by 
$(5 D_k, 5 H_k, S_k,V_k)$, and each $f \in \CF_k$ by $f^5$, one may assume
that 
\[ \sup_{|I| \leq S_k-1} |\log |\partial^I f|| \leq - 5 V_k \] 
for each $k$ sufficiently big and $f \in \CF_k$.

\satz{Definition}
Given the series $(D_k,H_k,S_k,V_k)$, and a $t \leq s-1$, an irreducible
subvariety $\CY$ of $\CX$ of codimension $p \leq t$
is called sufficiently approximating of order $k$ and multiplicity $S_Y \in \N$
with respect to $\theta \in \CX(\C)$, if the estimates
\begin{equation} \label{pcl2}
t_{H_k/D_k}(\CY) \leq \frac{S_Y}{S_k^p} 4^p D_k^{p-1} H_k, 
\end{equation}
and 
\begin{equation} \label{pcl3}
   \varphi_{H_k/D_k}^{S_Y/9^p}(\theta,\CY) \leq 
   -\frac{4 S_Y V_k}{14^{p-1} t_{\frac{H_k}{D_k}}(\CY) S_k}
\end{equation}
hold.
\end{Satz}

\satz{Lemma} \label{pcl}
Given the series $(D_k,H_k,S_k,V_k)$, let $C >> 0$ and $t \leq s-1$. 
Because of 
\[ \limsup_{k\to \infty} \frac{S_k^s V_k}{D_K^s(D_k+H_k)} = \infty \]
for every $k_0 \in \N$ there is a $k \geq k_0$ such that
\begin{equation} \label{C}
\frac{S_k^s V_k}{D_k^s(D_k + H_k)} \geq C,
\end{equation}
and assume $l \leq k$.
\begin{enumerate}

\item
Let $\CY$ be an irreducible subvariety of codimension $p$ in $\CX$, 
and $S_Y \in \N$ a number such that (\ref{pcl2}) holds.
Let further $f \in \CF_l$ be such that $\di f$ intersects $\CY$
properly, and assume
\[ D^{{S_Y} (S_l-1)/9^{p+1}}(\di f . Y,\theta) \leq 
   -\frac{4S_l S_Y V_k}{14^{p-1} S_k}. \]
Then, if either $k=l$ or $|\di f,\theta| \leq |Y,\theta|$,
there exists an irreducible component $\bar{\CY}$ of $\di f . \CY$ and
a number $S_{\bar{Y}}$ such that 
$S_{\bar{Y}}/t_a(\bar{\CY}) \geq S_Y/t_a(\di f . \CY)$, and
$\bar{\CY}$ is sufficiently approximating of order $k$ and multiplicity 
$S_{\bar{Y}}$ with respect to $\theta$.

\item
Let $\CY$ be an irreducible subvariety of codimension $p$ in $\CX$ that
is sufficiently approximating of order $k$ and multiplicity $S_Y$ with
respcet to $\theta$, and $f \in \CF_k$ a global section
whose restriction to $\CY$ is nonzero.
Then, there exists an irreducible component $\bar{\CY}$ of $\di f . \CY$,
and a number $S_{\bar{Y}} \in \N$ such that
$S_{\bar{Y}}/t_a(\bar{\CY}) \geq S_Y/t_a(\di f . \CY)$, and
$\bar{\CY}$ is sufficiently approximating of order $k$ and multiplicity 
$S_{\bar{Y}}$ with respect to $\theta$.

\end{enumerate}
\end{Satz}

\proof
1. Since 
\[ \varphi_{H_k/D_k}^{S_Y(S_l-1)/9^{p+1}} (\theta,\di f . \CY) \leq
   - \frac{4S_l S_Y V_k}{14^{p-1} t_{H_k/D_k}(\di f . \CY) S_k}, \]
Lemma \ref{gaal} implies that there is an irreducible component
$\bar{\CY}$ of $\di f . \CY$, and a number $S_{\bar{Y}}$ such that,
\[ \varphi_{H_k/D_k}^{S_{\bar{Y}}}(\bar{\CY},\theta) \leq 
   \varphi_{H_k/D_k}^{S_Y (S_l-1)}(f . \CY,\theta)
   + O(\log (D_k \deg Y)) \leq 
   -\frac{4\cdot S_l S_Y V_k}{4\cdot 14^{p-1}t_{H_k/D_k} (\di f . \CY)S_k}, \]
and by shrinking $S_{\bar{Y}}$ if necessary,
\begin{equation} \label{SY}
2 \; \frac{S_Y S_l}{t_{H_k/D_k}(\di f . \CY)} \geq
   S_{\bar{Y}}/t_{H_k/D_k}(\bar{\CY}) \geq 
   \frac{S_Y S_l}{t_{H_k/D_k} ( \di f . \CY)}. 
\end{equation}
Thereby,
\begin{equation}\label{varphi}
\varphi_{H_k/D_k}^{S_{\bar{Y}}} (\bar{Y},\theta) \leq 
   - \frac{4 S_{\bar{Y}} V_k}{14^p t_{H_k/D_k}(\CY)S_k}. 
\end{equation}

Further, by the algebraic and arithmetic B\'ezout Theorems, the inequality
$D_l < H_k$, and the fact that 
$\CY$ fulfills (\ref{pcl2}),
\begin{eqnarray*}
 t_{H_k/D_k}(\di f . \CY) &\leq&  D_l h(\CY) + H_l \deg Y + 
   \left(\frac{H_k}{D_k} +c\right) D_l \deg Y \\
&\leq& 2 D_l t_{H_k/D_k}(\CY) + 2 H_l \frac{D_k}{H_k} t_{H_k/D_k}(\CY) \\
&\leq& \frac{2S_Y}{S_k^p} 4^p D_l D_k^{p-1} H_k +
   \frac{2S_Y}{S_k^p} \frac{D_k}{H_k} 4^p H_l D_k^{p-1} H_k. 
\end{eqnarray*}
Hence, by the right hand side inequality of (\ref{SY}),
\[ \frac{S_{\bar{Y}}}{t_{H_k/D_k}(\bar{\CY})} \geq 
   \frac{S_l S_k^p}{2 \cdot 4^p D_l D_k^{p-1} H_k + 2 \cdot 4^p D_k^p H_l} \geq
   \frac{S_k^{p+1}}{4^{p+1} D_k^p H_k}, \]
the last inequality, because $l \leq k$ and both $D_k/S_k$ and $H_k/D_k$
are non-decreasing. Thereby,
\[ t_{H_k/D_k} (\CY) \leq \frac{S_{\bar{Y}}}{S_k^{p+1}} 4^{p+1} D_k^p H_k, \]
that is $\bar{\CY}$ is sufficiently approximating of order $k$ and multiplicity 
$S_{\bar{Y}}$ with respect to $\theta$.

\vspace{2mm}

2. For $k=l$, since $\di f$ intersects $\CY$ properly, 
by the derivative metric B\'ezout Theorem (\ref{DMBT}),
\begin{eqnarray*} 2D^{S_Y (S_k-1)/9^{p+1}}(\di f. Y,\theta) &\leq& 
   \mbox{max}(S_kD^{S_Y/9^p}(Y,\theta),S_YD^{(S_k-1)/9^p}(\di f, \theta)) \\ 
&+& 2H_k \deg Y + 2D_k h(\CY) + 2d D_k \deg Y \\
&+& c (D_k \deg Y) \log (D_k \deg Y) \\ \\
&\leq& \mbox{max}(S_kD^{S_Y/9^p}(Y,\theta),S_YD^{(S_k-1)/9^p}(\di f, \theta)) \\
&+&  7 D_k t_{H_k/D_k} (\CY) \log(D_k \deg Y). 
\end{eqnarray*}
Further, by (\ref{pcl3}), and Proposition \ref{glschn},
\[ S_k D^{S_Y/9^p}(Y,\theta) \leq - \frac{4 S_Y V_k}{14^{p-1}}, \]
\[ S_Y D^{(S_k-1)/9^p} (\di f, \theta) \leq - 5 S_Y V_k + c D_k \log D_k \leq
   - 4 S_Y V_k, \]
and by (\ref{pcl2}) and (\ref{C}), since $p \leq t \leq s-1$,
\begin{eqnarray*} 
  7 D_k t_{H_k/D_k} (\CY) \log (D_k \deg Y) &\leq& 
   7 \cdot 4^p \frac{S_Y}{S_k^p} D_k^p H_k \log (D_k \deg Y) \\
&\leq& 7 \cdot 4^p S_Y V_K/C \log (D_k \deg Y) \leq
   \frac{S_Y V_k}{14^{p-1}}, 
\end{eqnarray*}
for $C$ sufficiently big.
Hence,
\[ 2D^{S_Y (S_k-1)/9^{p+1}}(\di f. Y,\theta) \leq 
   -\frac{4 S_Y V_k}{2 \cdot 14^{p-1}} = 
   - \frac{4 \cdot S_Y S_k V_k}{2 \cdot 14^p S_k}, \]
that is
\[ \varphi_{\frac{H_k}{D_k}}^{S_Y (S_k-1)/9^{p+1}}(\di f . \CY,\theta) \leq 
   -\frac{4 S_k S_Y V_k}{2 \cdot 14^p S_k t_{\frac{H_k}{D_k}}}(\di f . \CY). \]
Thereby the premisses of part 1 are fulfilled with $l=k$, and part one implies 
the claim.

If $l<k$, and $|\di f,\theta| \leq |Y,\theta|$, the claim follows similarly, this
time using Corollary \ref{DMBTcor}.2.

\proof {\sc of Theorem \ref{algind1}}
Assume $t\leq s+1$, let $k_0 \in \N$ be any number, and
\[ R = \inf \{ \log |\di f,\theta| \; | f \in \Gamma(\Pe^t,O(D_{k_0})),
        \log |f| \leq H_{k_0},  f  \neq 0 \}. \]
Let further $C$ be an arbitrarily big constant, and $k > k_0$ such that
\begin{equation} \label{gross} 
\frac{S_k^s V_k}{D_k^s (D_k + H_k)} \geq C, 
\end{equation}
and
\[ \frac{S_k^{t-1} V_k}{D_k^{t-1} (D_k + H_k)} \geq C R. \]
Let $\CY \subset \CX$ be a subvariety of maximal codimension that is
sufficiently approximating of order $k$ and some multiplicity $S_Y$.
Then $Y$ is contained in the support of $\di f$ for every $f \in \CF_k$, since
otherwise, by Lemma \ref{pcl}.2, there would be a subvariety $\bar{\CY}$ of
codimension $p+1$ fulfilling the same conditions, thereby contradicting
the maximality of the codimension of $\CY$. 
Since the derivated algebraic distance of the zero cycle is defined as $0$,
we have $p \leq t$.

Let now
\[ l = \mbox{max} 
   \{ \bar{k} \leq k | \exists f \in \CF_{\bar{k}-1} : Y \not\subset
   supp(\di f) \}, \]
Then, $Y$ is contained in the support of $\di f$ for
every $f \in \CF_l$, hence
\begin{equation} \label{vunten}
\log |Y,\theta| > -V_{l-1}/(S_{l-1}), 
\end{equation}
and for every $f \in \CF_l$, by \cite{App1}, Theorem 2.2.1 and 
(\ref{pcl2}), and (\ref{pcl3}),
\[ \log |\di f,\theta| \leq \log |Y,\theta| \leq 
   \varphi_{H_k/D_k}(Y,\theta) + c \leq
   \varphi_{H_k/D_k}^{S_Y/9^p}(Y,\theta) \leq \]
\[ - \frac{4 S_Y V_k}{14^{p-1}t_{H_k/D_k}(\CY) S_k} \leq
   - \frac{4 S_Y V_k S_k^{p-1}}{14^{p-1}S_k S_Y D_k^{p-1}} H_k \leq
   - \frac{4 V_k S_k^{t-1}}{14^{p-1}  D_k^{-t}(D_k+H_k)}
     \frac{S_{k_0}^{t-p}}{D_{k_0}^{t-p}} < -R, \]
the last inequality holding if the constant $C$ is chosen sufficiently big.
The inequalities $\log |di f, \theta| < -R$ for every $f \in \CF_l$ imply
$l > k_0$.

Let $D = [S_{l-1} S_Y V_k/(14^{p-1} V_{l-1})]$. If
$\deg Y \leq D/3$, then, again by \cite{App1}, Theorem 2.2.1,
\[ \log |Y,\theta| \leq \frac{- S_Y V_k}{3 \cdot14^{p-1} D} \leq 
   - V_{l-1}/3S_{l-1}\]
in contradiction with (\ref{vunten}). If $\deg Y \geq D$, let 
$g \in \CF_{l-1}$ be such that $Y \not\subset \mbox{supp} (\di g)$. 
If $|\di g,\theta| \leq |Y,\theta|$, Lemma \ref{pcl}.1 would contradict
the minimality of the dimension of $\CY$. 
Hence, $|Y, \theta| \leq |\di g, \theta|$.

Using Corollary \ref{DMBTcor} for 
$Z_0=Y, Z_1 = \di g, d_0 = D, S_0 =S_Y, S_1 = S_{l-1}$, one gets a 
$K \leq D S_{l-1}$ such that
\[ K \log |Y, \theta| + 2D^{S_Y (S_{l-1}-1)/9^{p+1}} (Y. \deg g, \theta) \leq\]
\[ \mbox{max}(D\;  D^{S_{l-1}-1}(\di g, \theta), S_{l-1} D^{S_Y}(Y,\theta) + \]
\[ 2 H_{l-1} \deg Y + 2 D_{l-1} h(\CY) + 2 d D_{l-1} \deg Y. \]
Since 
$D^{S_{l-1}-1}(\di g,\theta) \leq -5V_{l-1}, 
D^{S_Y}(Y,\theta) \leq -4S_YV_k/14^{p-1} S_k$,
and by assumption $H_{l-1}/S_{l-1} \leq H_k/S_k$, and 
$D_{l-1}/S_{l-1} \leq D_k/S_k$, the
above is less or equal
\[ \mbox{max}(-5S_{l-1} 
   S_Y V_k/(2 \cdot 14^{p-1}), -S_{l-1} S_Y V_k/( 14^{p-1}))+ \]
\[ H_k \frac{S_{l-1}}{S_k} \deg Y +
   D_k \frac{S_{l-1}}{S_k} h(\CY) + d D_k \frac{S_{l-1}}{S_k} \deg Y. \]
Further, by (\ref{pcl2})
\begin{eqnarray*}
   2D_k \frac{S_{l-1}}{S_k} h(\CY) \leq 2D_k \frac{S_{l-1}}{S_k} 
   t_{\frac{H_k}{D_k}}(\CY) &\leq&
   2D_k \frac{S_{l-1}}{S_k} \frac{4^p S_Y D_k^{p-1}(D_k+H_k)}{S_k^p} \\
&=& 2S_{l-1} S_Y \frac{4^p D_k^p (D_k+H_k)}{S_k^{p+1}} \leq 
   2 \cdot 4^p S_{l-1} S_Y \frac{V_k}C.
\end{eqnarray*}
The last inequality because of $p \leq t$.
Similarly,
\[ 2H_k \frac{S_{l-1}}{S_k} \deg Y \leq 
   2 \cdot 4^p S_{l-1} S_Y \frac{V_k}C, \quad 
   d D_k \frac{S_{l-1}}{S_k} \deg Y \leq 2 \cdot 4^p S_{l-1} S_Y \frac{V_k}C.\]
Hence,
\[ K \log |\di g + Y, \theta| + 2D^{9S_Y S_{l-1}/9} (Y. \deg g, \theta) \leq \]
\[ - 5 S_{l-1} S_Y V_k/ (2 \cdot 14^{p-1}) + 6 S_{l-1} S_Y V_k/C \leq
   - S_{l-1} S_Y V_k/ (2 \cdot 14^{p-1}), \]
for $C$ sufficiently large.

Since $\CY$ was chosen of maximal codimension, Lemma \ref{pcl}.1 implies \\
$D^{S_Y S_{l-1}/9} (Y. \di g ,\theta) \geq 
- S_{l-1} S_Y V_k/ (4 \cdot 14^{p-1})$.
Consequently,
\[ K \log |\di g + Y, \theta| \leq - S_{l-1} S_Y V_k/ (4 \cdot 14^{p-1}), \]
and thereby
\[ \log |Y,\theta| \leq - S_{l-1} S_Y V_k/ (4K \cdot 14^{p-1}). \]
Since $K \leq S_{l-1} D$, this is less or equal
\[ -S_Y V_k/(4D 14^{p-1}) \leq -V_{l-1}/(4 S_{l-1}), \]
again contradicting (\ref{vunten}). Since the assumtions
$t-1 \leq s$ leads to a contradiction, we have $t-1 > s$.

\vspace{2mm}

\section{Proof of second criterion}

To prove Theorem \ref{algind2}, let $\theta$ be a point in projective
space $\Pe^M$, assume its algebraic closure $\CX$ over $\spec \; \Z$ has
relative dimension $t$, and let $(D_k,S_k,H_k,V_k)$ be a quadrupel
of series fulfilling the assumptions of the Theorem. Let further
$F,G$ be the functions $F(k)=D_k/S_k, G(k) = H_k/D_k$. 
Since $F,G$ are of uniform polynomial growth, by Lemma \ref{upg}, there
is a $k_0$ such that for every $k \geq k_0$,
\begin{equation} \label{rega}
\frac12 D_{k+1}/S_{k+1} \leq D_k/S_k \leq D_{k+1}/S_{k+1}, \quad
\frac12 H_{k+1}/S_{k+1} \leq H_k/S_k \leq H_{k+1}/S_{k+1}. 
\end{equation}
By Lemma \ref{upg}, the function $H(D) = G \circ F^{-1}(D)$ is of uniform 
polynomial growth with $n_H \geq 0$. Multiplying $H(D)$ by a positive 
constant, if necessary, one can assure that $H(D) \geq a D$ with an arbitrary 
number $a \geq 1$.
By Proposition \ref{2hilf}, there are numbers $b_1,1>c_0>0,n_1 \in \N$ 
and an infinite subset $M \subset \N$ such that for each $D \in M$ there is
an irreducible variety $\beta_{nD}$ of codimension $t$ in $\Pe^t$ and
a locally complete intersection $\CZ$ at $\alpha_{nD}$ of codimension 
$r \leq t-1$, such that
\[ \deg \beta_{n_1D} \leq (n_1D)^t, \quad h(\beta_{n_1D}) \leq H(n_1D) (n_1D)^{t-1}, 
   \quad D(\beta_{n_1D}, \theta) \leq - b_1 t_H(\beta_{n_1D}) D, \]
\begin{equation} \label{rega-1}
t_{H/D}(\beta_{n_1D}) \geq c_0 t_{H/D}(\CZ_{min}) D^{t-r},
\end{equation}
where $\CZ_{min}$ is the irreducible component of  $\CZ$ with
minimal $\frac HD$-size.
Let $\pi_X \to \Pe^t$ be the projection from section \ref{psp}, and
$\alpha_D \subset \CX$ an irreducible component of $\pi_X^* \alpha_D$,
further $\CY$ an irreducible component of $\pi_X^* \CZ_{min}$ containing
$\alpha_D$. By (\ref{rega-1}), Proposition \ref{hoch}.1, and 
Proposition \ref{runter}, there are constants
$b, 1 > c >0, n \in \N$ such that
\[ \deg \alpha_{nD} \leq (nD)^t, \quad h(\alpha_{nD}) \leq H(nD) (nD)^{t-1}, 
   \quad D(\alpha_{nD}, \theta) \leq - b t_H(\alpha_{nD}) D, \]
\begin{equation} \label{rega0}
t_{H/D}(\alpha_{nD}) \geq c t_{H/D}(\CY) D^{t-r}.
\end{equation}
With a big constant $c_3$ put
\[c_1 = \frac{c}{9M (h(\CX) + c_3 \deg X)}. \]
Since $\lim_{k \to \infty} \frac{V_kS_k^s}{D_k^s(D_k+H_k)} = \infty$, there is a
$k_1 \geq k_0$ such that
\[ \frac{V_k S_k^s}{D_k^s(D_k + H_k)} > 
   40 M h(\CX+ c_3 \deg X)(d+1) (2 n \; \mbox{max}(1/c_1,(10+d)/b))^t.\]
for every $k \geq k_1$, where $d$ is the constant from Proposition 
\ref{liou}. Since $M$ is infinite, (\ref{rega})
implies that there is a $D \in M$ and a $k \geq k_1$ such that 
\begin{equation} \label{rega1}
\left( \frac{\mbox{min}(c_1,b/(10+d))}2\right) D \leq \frac{D_k}{S_k} 
   < \left( \mbox{min}(c_1,b/(10+d)) \right) D. 
\end{equation}
Applying the function $H = G \circ F^{-1}$ to both sides, and using that
it is eventually non-decreasing, gives
\begin{equation} \label{rega2}
\left( \frac{\mbox{min}(c_1,b/(10+d))}{2} \right) H \leq 
   \frac{H_k}{S_k} <
   \left( \mbox{min}(c_1,b/(10+d)) \right) H, 
\end{equation}
with $H = H(D)$. Adding both inequalities implies
\[ \left( \frac{\mbox{min}(c_1,b/(10+d))}{2} \right) (H + D) \leq 
   \frac{H_k+D_k}{S_k} < \mbox{min}(c_1,b/(10+d))  (H+D) \leq \]
\begin{equation} \label{rega3} 
2 \mbox{min} (c_1,b/(10+d))H. 
\end{equation}

For a given global section $h \in \Gamma(\Pe^M,O(1))$ with $h_\theta \neq 0$,
identify an $f \in \CF_k$ with $f/h^{D_k} \in \Q(X)$.

\satz{Lemma}
There is an $f \in \CF_k$ such that for 
some $I$ with $|I| \leq 2 S_k/3$ the restriction of $\partial^I f$
to $\alpha_{nD}$ is nonzero.
\end{Satz}

\proof
Assume the opposite, and inductively construe
a chain of subvarieties
\[ \CY_1 \supset \cdots \supset \CY_{t-r} = \alpha_{nD}, \]
such that
\[ \frac{t_{\frac HD}(\CY_i)}{v_{\alpha_{nD}}(Y_i)} \leq c^{i-1} D^{i-1} 
   t_{\frac HD}(\CY), \quad i=1, \ldots {t-r}, \]
\[ \frac{t_{\frac HD}(\CY_i)}{v_{\alpha_{nD}}(Y_i)} \leq c^{i-1} D^{i-1} 
   t_{\frac HD}(\CY), \quad i=2, \ldots {t-r}, \]
in the following way: Since $\alpha_{nD}$ is contained in 
$Y$, by fact \ref{dim0}, we have $v_{\alpha_{nD}}(Y) \geq 1$, thus can
choose $\CY_1 = \CY$. Assume $\CY_i$ is given, and fulfills the above estimate.
Since $\alpha_{nD}$ is contained in $\CY_i$,
by \cite{App1}, Theorem 2.2.2,
\begin{equation} \label{nah}
 \log|\CY_j,\theta| \leq \log |\alpha_{nD},\theta| \leq
   \frac{D(\alpha_{nD},\theta)}{t_{\frac HD}(\alpha_{nD})} + O(1) \leq
   - b D +O(1). 
\end{equation}
Thus, for $k$ sufficiently large, the assumption in the Theorem asserts that 
there is an $f_i \in \CF_k$ and a multi index $I_i$ with $|I_i| \leq S_k/3$
such that the restriction of 
$\partial^{I_i} f_i$ to $\CY_i$ is nonzero, and by Theorem \ref{Abl}, there
are polynomials $P, f_{I_i}$ with $P|_X \neq 0$, thereby $P(\theta) \neq 0$,
and by (\ref{nah}) also $P(\alpha_{nD}) \neq 0$, fulfilling 
\[ \deg f_{I_i} \leq \deg f_i + (2S-1)(M-t) \deg X \leq
   2 M \deg X D_k, \]
\begin{eqnarray*} 
\log |f_{I_i}| &\leq& \log |f_i| + \log \deg f_i \\
&+&   (2S-1)(M-t) (h(\CX) + c_4 \deg X +\log \deg X) + \log (2S!) \\ \\
&\leq& 
(2 M (h(\CX)+ c_3 \deg X) (H_k + D_k) \leq 3 M (h(\CX)+ c_3 \deg X) H_k, 
\end{eqnarray*}
and
\[ \partial^{I_i} f= \frac{f_{I_i}}{P^{2|I_i|-1}}, \]
and thereby
\[ \partial^J f_{I_i}(\alpha_{nD}) = 
   \partial^J (\partial^{I_i} f P^{2|I_i|-1})(\alpha_{nD}) = 0 \]
for every $J$ with $|J| \leq S_k/3$.
Hence, by Proposition \ref{mult} $v_{\alpha_{nD}}(\di f_{I_i}) \geq S_k/3$, 
and by the local B\'ezout Theorem,
\[ v_{\alpha_{nD}}(Y_i . \di f_{I_i}) \geq \frac{S_k}3 v_{\alpha_{nD}}(Y_i). \]
Further, by the algebraic B\'ezout Theorem,
\[ \deg (Y_i . \di f_{I_i})\leq 2 M \deg X D_k \deg Y_i, \]
\[ h(\CY_i. \di f_{I_i}) \leq \]
\[ 2 M \deg X D_k h(\CY_i) + 
   3 M (h(\CX) + c_3 \deg X) H_k \deg Y_i +
   2 c M \deg X D_k \deg Y_i \leq \]
\[ 2 M \deg X D_k h(\CY_i) + 
   4 M (h(\CX) + c_3 \deg X) H_k \deg Y_i. \]
Hence,
\[ t_{\frac HD}(\CY_i .\di f_{I_i}) \leq \]
\[ 2 \frac HD M \deg X D_k \deg Y_i +
   2 M \deg X D_k h(\CY_i) + 
   4 M (h(\CX) + c_3 \deg X) H_k \deg Y_i \leq \]
\[ 2 M \deg X D_k t_{\frac HD} (\CY_i) +
   4 M (h(\CX) + c_3 \deg X) \frac{H_k}{H} D \; t_{\frac HD} (\CY_i). \]
Together with the above estimate on the order of vanishing of 
$Y_i . \di f_{I_i}$ at $\alpha_{nD}$, this gives
\[ \frac{t_{\frac HD}(\CY_i. \di f_{I_i})}{v_{\alpha_{nD}}(Y_i .\di f_{I_i})} \leq
   \frac3{S_k} \left(2 M \deg X D_k + 4 M (h(\CX)+c_3 \deg X) \frac{H_k}HD\right)
   \frac{t_{\frac HD}(\CY_i)}{v_{\alpha_{nD}}(Y_i)}, \]
which by (\ref{rega1}), and (\ref{rega2}) is less or equal
\[ \left(\frac{\mbox{min}(c_1,b)}2\right) D \left( 2 M \deg X +
   4 M (h(\CX) + c_3 \deg X) \right)
   \frac{t_{\frac HD}(\CY_i)}{v_{\alpha_{nD}}(Y_i)} < \]
\[ \left(\frac{\mbox{min}(c_1,b)}2\right) D \left( 3 M \deg X +
   5 M (h(\CX) + c_3 \deg X) \right)
   \frac{t_{\frac HD}(\CY_i)}{v_{\alpha_{nD}}(Y_i)} \leq \]
\[ c D \frac{t_{\frac HD}(\CY_i)}{v_{\alpha_{nD}}(Y_i)}. \]
By Lemma \ref{wov}, there is an irreducible component $\CY_{i+1}$
of $\CY_i \di f_{I_i}$ such that
\[ \frac{t_{\frac HD}(\CY_{i+1})}{v_{\alpha_{nd}}(Y_{i+1})} \leq
   \frac{t_{\frac HD}(\CY_i . \di f_{I_i})}{v_{\alpha_{nd}}(Y_i. \di f_{I_i})}, \]
which by the above is less than
\[ c D \frac{t_{\frac HD}(\CY_i)}{v_{\alpha_{nD}}(Y_i)}, \]
which by induction hypothesis is less or equal
\[ c^i D^i t_{\frac HD}(\CY). \]
proving the claim for $i+1$.
For $i = t-r$, the claim gives
\[ t_{\frac HD}(\alpha_{nD}) = \frac{t_{\frac HD}(\alpha_{nD})}
   {v_{\alpha_{nD}}(\alpha_{nD})} < c^{t-r} D^{t-r} t_{\frac HD}(\CY), \]
contradicting the lower estimate on $t_{\frac HD}(\alpha_{nD})$ in
(\ref{rega0})

\proof {\sc of Theorem \ref{algind2}, continuation}
Let $g = \partial^I f$ with $|I| \leq 2S_k/3$ be as in the Lemma. By
Theorem \ref{Abl}, there are polynomials $P,g_I$ such that with $c_3$
chosen sufficiently big,
\[ g = \frac{g_I}{P^{2|I|-1}}, \quad \deg g_I \leq 2 M \deg X D_k, \]
\[ \log|g_I| \leq 3 M (h(\CX) + c_3 \deg X) H_k, \]
and by Corollary \ref{distAbl},
\[ \sup_{|J| \leq S_k/3} \log |\la \partial^I g_I|\theta\ra| \leq -V_k/2. \] 
Theorem \ref{Abl} also implies $P|_X\neq 0$, and thereby $P(\theta) \neq 0$,
which by (\ref{nah}) implies 
\[ P(\alpha_{nD}) \neq 0. \]
Hence, by Theorem \ref{DMBT}.2, and Corollary \ref{distAbl},
\[ D(\alpha_{nD},\di g_I) \leq
   \mbox{max}\left(\frac{S_k}3 D(\alpha_{nD},\theta),-V_k/2\right) \leq
   \mbox{max}\left(-b \frac{S_k}
   3 D t_{H_D}(\alpha_{nD}),-V_k/2\right). \]
Further, by Liouvilles Theorem \ref{liou},
\begin{eqnarray*} D(\alpha_{nD},\di g_I) &\geq& - 2 M \deg X D_k h(\alpha_{nD}) - 
   \deg \alpha_{nD}  3 M(h(\CX) + c_3 \deg X) H_k \\ 
&& -   2 d M  \deg X D_k \deg \alpha_{nD} \\ \\
&\geq&- 2 M \deg X D_k t_{\frac HD} (\alpha_{nD} ) - 
   3 M (h(\CX) + c_3 \deg X) \frac{H_k D}H t_{\frac HD} (\alpha_{nD}) \\
&& -2 d M \deg X D_k \frac DH t_{\frac HD} (\alpha_{nD}) \\ \\
&\geq& - \left(2(d+1)M \deg X D_k + 3 M (h(\CX) + c_3 \deg X)
   \frac{H_k D}H\right) \times \\ && t_{\frac HD}(\alpha_{nD}). 
\end{eqnarray*}
The two inequalities together give
\[ - \left(2(d+1)M \deg X D_k + 3 M (h(\CX) + c_3 \deg X)\frac{H_k D}H\right) 
   t_{\frac HD}(\alpha_{nD}) \leq \]
\[ \mbox{max}\left(-b 
   \frac{S_k}3 D t_{\frac HD}(\alpha_{nD}),-V_k/2\right). \]
If $-(2 M (d+1)\deg X D_k + 3 M (h(\CX) + c_3 \deg X)\frac{H_k D}H) 
t_{H_D}(\alpha_{nD})$ were less or equal
$-b(S_k/3) D t_{\frac HD}(\alpha_{nD})$, if $c_3$ is chosen sufficiently big,
this would contradict the second inequality of (\ref{rega3}). 
Hence,
\begin{equation} \label{schl}
-\left(2(d+1)M \deg X D_k + 3 M (h(\CX)+c_3 \deg X)
\frac{H_k D}H\right) t_{H_D}(\alpha_{nD}) \leq -V_k/2. 
\end{equation}
By the upper estimates on $\deg \alpha_{nD}$, and $h(\alpha_{nD})$,
\[ \left(2(d+1)M \deg X D_k + 3 M (h(\CX) + c_3 \deg X)\frac{H_k D}H\right) 
   t_{H_D}(\alpha_{nD}) \leq \]
\[ 2( 2 M \deg X (d+1)D_k + 3 M (h(\CX + c_3 \deg X)
   \frac{H_k D}H) 2H n^t D^{t-1}, \]
which by (\ref{rega1}) and (\ref{rega2}) is less or equal
\[ 8 M \deg X (d+1) (2n)^t \mbox{max}(1/c_1,(10+d)/b))^t\frac{H_k D_k^t}{S_k^t}+\]
\[ 12 M (h(\CX + c_3 \deg X)(2n)^t 
   \mbox{max}(1/c_1,(10+d)/b))^t \frac{H_k D_k^t}{S_k^t} \leq \]
\[ 18 M h(\CX + c_3 \deg X) (d+1)(2n^t) \mbox{max}(1/c_1,(10+d)/b))^t 
   \frac{H_k D_k^t}{S_k^t}, \]
for $c_3$ sufficiently big.
Together with (\ref{schl}), this implies
\[ \frac{V_k S_k^t}{D_k^t H_k} \leq 
   40 M h(\CX + c_3 \deg X)(d+1) (2 n \; \mbox{max}(1/c_1,(10+d)/b))^t. \]
Since $k$ was chosen such that
\[ \frac{V_k S_k^s}{D_k^s(S_k+H_k)} > 40 M (h(\CX)+c_3 \deg X) (d+1) 
   (2 n \; \mbox{max}(1/c_1,(10+d)/b))^t, \]
and $S_k/D_k <1$, we get $t \geq s+1$.


\begin{thebibliography}{0mm}

\bibitem[Be]{Be} D.\@ Bertrand: Duality on Tori and multiplicative
dependence relations. J.\@ Austral.\@ Math.\@ Soc.\@ (Series A) 
62 (1997), 198--216.

\bibitem[BGS]{BGS} Bost, Gillet, Soul\'e: Heights of projective varieties and
positive Green forms. JAMS 7,4 (1994)

%\bibitem[Ch]{Ch} M.\@ Chardin: Une majoration de la fonction de Hilbert et
%ses cons\'equences pour l'interpolation alg\'ebrique. Bull.\@ Soc.\@
%Math.\@ France 117 (1989) 305--318.

%\bibitem[CP]{CP} M.\@ Chardin, P.\@ Philippon: R\'egularite et interpolation.
%J.\@ of alg.\@ geom.\@ 8 (1999), 471-481.


\bibitem[LR]{LR} M.\@ Laurent, D.\@ Roy: Criteria of algebraic independence
with multiplicities and approximation by hypersurfaces.

\bibitem[Ma1]{App1} H.\@ Massold: Diophantine Approximation 
on varieties I: Algebraic Distance and Metric B\'ezout Theorem. 
math.NT/0611715v2

\bibitem [Ma2]{App2} H.\@ Massold: Diophantine Approximation
on varieties II: Explicit estimates for arithmetic Hilbert functions. 
0711.1667


\bibitem[Ma3]{App3} H.\@ Massold: Diophantine Approximation on varieties III;
   Approximation of non algebraic points by algebraic points. arxiv: 
   0711.3645

\bibitem[Ma4]{App4} H.\@ Massold: Diophantine Approximation on varieties IV:
 Derivated algebraic distance and derivative metric  B\'ezout Theorem. 
  arxiv: 0901.3889

\bibitem[Ma5]{Mahler} H.\@ Massold: Mahler classification for points on
algebraic varieties. To appear

%\bibitem[Ma6]{trans} H.\@Massold:    To appear.

\bibitem[Ma6]{Liouville} H.\@ Massold: Liouville inequality for arbitrarily
dimensional subvarieties. To appear.


\bibitem[Ph]{Ph} P.\@ Philippon, Approximations alg\`ebrique dans les
espaces projectifs I. Journal of Number Theory 81, (2000) 234--253.

%\bibitem[R]{R} D.\@ Roy: Approximation alg\`ebrique simultan\'ee de nombres
%    de Liouville. 


\bibitem[RW]{RW} D.\@ Roy, M.\@ Waldschmidt: Approximation diophantienne
et ind\'ependance alg\'ebrique de logarithmes. Ann.\@ sc.\@ de l'ENS
30,6 (1997)

\bibitem[SABK]{SABK} C.\@ Soul\'e, D.\@ Abramovich, J.\@-F.\@ Burnol,
               J.\@ Kramer: Lectures on Arakelov Geometry. 
               Cambridge University Press 1992 

\end{thebibliography}
\end{document}